\documentclass[twocolumn]{autart}    

\usepackage{graphicx}          
\usepackage{amssymb, amsmath, cite}
\usepackage{array}
\usepackage{enumerate}
\usepackage{tikz}
\usepackage{caption}
\usepackage{float}
\usetikzlibrary{arrows}
\usepackage[Symbol]{upgreek}
\newtheorem{theorem}{Theorem}
\newtheorem{definition}{Definition}
\newtheorem{example}{Example}
\newtheorem{lemma}{Lemma}
\newtheorem{remark}{Remark}
\newcommand{\CommaBin}{\mathbin{\raisebox{0.5ex}{,}}}

\newcommand\ddfrac[2]{\frac{\displaystyle #1}{\displaystyle #2}}
  {
     \theoremstyle{definition}
      \newtheorem{assumption}{Assumption}
  }
\newcommand{\RomanNumeralCaps}[1]
    {\MakeUppercase{\romannumeral #1}}
\begin{document}

\begin{frontmatter}

\title{Negative Imaginary Systems Theory for Nonlinear Systems: A Dissipativity Approach \thanksref{footnoteinfo}} 

\thanks[footnoteinfo]{This paper was not presented at any IFAC
meeting. Corresponding author Ahmed~G.~Ghallab.} 
\thanks{This work was supported by the Australian Research Council under grant DP190102158.}

\author[ACT]{Ahmed~G.~Ghallab}\ead{ahmed.ibrahim@anu.edu.au} ~~and~~   
\author[ACT]{Ian~R.~Petersen}\ead{ian.petersen@anu.edu.au}  

\address[ACT]{Research School of Engineering, The Australian National University, Canberra 0200, Australia}  


\begin{abstract}                          
Negative imaginary (NI) systems theory is a well-established system theoretic framework for analysis and design of linear-time-invariant (LTI) control systems. In this paper, we aim to generalize negative imaginary systems theory to a class of nonlinear systems. Based on the time domain interpretation of the NI property for LTI systems, a formal definition in terms of a dissipativity with an appropriate work rate will be used to define the nonlinear negative imaginary (NNI) property for a general nonlinear system. Mechanical systems with force actuators and position sensors are nonlinear negative imaginary according to this new definition. Using Lyapunov stability theory, we seek to establish a nonlinear generalization of the NI robust stability result for positive feedback interconnections of NNI systems. An example of a nonlinear mass-spring-damper system with force as input and displacement of the mass as output will be presented to illustrate the applicability of the NNI stability result. Furthermore, the case of NI systems with free motion will be investigated in the nonlinear domain based on the dissipativity framework of NNI systems.
\end{abstract}

\end{frontmatter}

\section{Introduction}
One of the most appealing tools in the field of nonlinear control design is the "classical" dissipativity and passivity theory which characterizes the dissipation of energy in the system with respect to a supplied energy rate from the environment (see \cite{willems1972dissipative1, willems1972dissipative2, brogliato-bk2007, van2000l2} for detailed accounts). However, many systems that dissipate energy in the physical sense don't fall into this classical framework. Concretely speaking, many real-world systems can be considered where the input $u$ is some form of mechanical or electrical force, and the output $y$ is a corresponding displacement. In this case, the associated energy pair is not $(u,y)$, but $(u,\dot{y})$ and the systems are dissipative (passive) from $u$ to $\dot{y}$. For instance, flexible structures with colocated force actuators and position sensors are passive (or dissipative) from the input to the derivative of the output (instead of the output as in the classical passivity theory). For further examples of systems which are passive from the input to the time-derivative of the output were investigated in \cite{Angeli2006, Padthe2005, gorbet1998generalized, gorbet2001passivity, morris1994dissipative, pare2001kyp, brogliato-bk2007}.

As can be seen from these examples, the analysis and robust control of such systems cannot be handled directly using the classical dissipativity and passivity theory where, on one hand, systems with rate sensors and force actuators often violate the standard sector condition for passivity. On the other hand, the restriction of collocated rate sensors and actuators is a stringent one, and excludes many applications. These observations motivates the need for an extended framework of the dissipativity and passivity of the systems to allow for more general supply rates which involve derivatives of the inputs and outputs. This in turn will enable replacing rate sensors on a structure with displacement sensors.


In this regard, negative imaginary systems theory has proven to be an effective tool in the analysis and control design of positive real LTI systems \cite{petersen2010, lanzon2017feedback}. NI systems theory extends and complements range of applicability of passivity theory for LTI systems with relative degree greater than one. An NI~LTI system is a passive system from the input to the derivative of the output. One of the fruitful results in the existing NI literature pertaining to the stability robustness of positive feedback interconnections of LTI systems, as shown in Figure \ref{sys_31}, is given in \cite{petersen2010}.
This LTI NI stability result shows that a necessary and sufficient condition for internal stability when the plant is NI, with transfer function matrix $H_1(s)$, and the controller is a strictly negative imaginary (SNI) system, with transfer function matrix $H_2(s)$, is given by the following DC-gain condition:
\begin{equation}\label{dcloop}
\lambda_{max}(H_1(0)H_2(0)) < 1,
\end{equation}
where $\lambda_{max}$ denotes the maximum eigenvalue of a matrix with only real eigenvalues. This means that the internal stability of such a feedback system is governed by the loop gain at steady-state.

This NI stability result is of practical importance and has been employed in a wide variety of control applications including, for instance, robust vibration control of flexible structures, atomic force microscopy, and nano-positioning systems. Also, NI systems theory has been used for the robust control of highly resonant flexible structures with colocated position sensors and force actuators \cite{petersen2010, lanzon2008, Mabrok2015b, CaiAug.2010, Bhikkaji2012, Mabrok2011a, Mabrok2015d}. For more exposure on applications of NI theory in the literature see, for instance, \cite{Bhikkaji2009, CaiAug.2010, Diaz2012, Mahmood2011}.

\tikzstyle{block} = [draw, thick, rectangle,
    minimum height=4em, minimum width=4em, samples=501]
\tikzstyle{sum} = [draw, circle,inner sep=0pt,minimum size=1pt, node distance=1cm]
\tikzstyle{input} = [coordinate]
\tikzstyle{output} = [coordinate]
\tikzstyle{pinstyle} = [pin edge={to-,thin,black}]
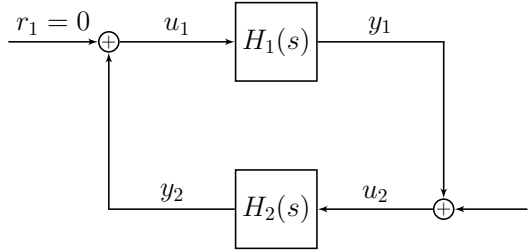
\begin{figure}[t]
\centering
\setlength{\belowcaptionskip}{-16pt}
\resizebox{7 cm}{!}{
\begin{tikzpicture}[auto, thick, node distance=2cm,>=latex']
    \node [input, name=input] {};
    \node [sum, right of=input, node distance=1.8cm] (sum1) {$+$};
    \node [block, right of=sum1,
            node distance=3cm] (system) {\large $H_1(s)$};
\node[inner sep=0,minimum size=0,right of=system,node distance=3cm] (k) {}; 
    \node [block, below of=system, node distance=3cm] (measurements) {\large $H_2(s)$};
    \node[inner sep=0,minimum size=0,left of=measurements,node distance=3cm] (p) {}; 
    \node [sum, below of=k, node distance=3cm] (sum2) {$+$};
    \node [input,  right of=sum2, node distance=1.5cm](input2) {};
    \draw [samples=501,draw,->] (input) -- node {\large $r_1=0$} (sum1);
    \draw [->] (sum1) -- node {\large $u_1$} (system);
    \draw [-] (system) -- node {\large $y_1$}(k);
    \draw [->] (k) -| node {} (sum2);
    \draw [->] (input2) -- node[above]{} (sum2); 
    \draw [->] (sum2) -- node[above]{\large $u_2$} (measurements);
    \draw [-] (measurements) -- node[above]{\large $y_2$} (p);
    \draw [->] (p) -|  node {} (sum1);
\end{tikzpicture}}
\caption{Feedback interconnection of two LTI systems where the plant is NI with transfer function $H_1(s)$ and the controller is SNI with transfer function $H_2(s)$}\label{sys_31}
\end{figure}

In real-world applications, most physical and engineering systems exhibit nonlinear behaviour, which in general makes a controller design difficult. We aim in this paper to develop an energy-based control framework for the analysis and design of dissipative (passive) system with supply rates involve the time derivative of input and output of the system. In the literature, different approaches have been proposed to deal with  analysis and control design of such type of systems, see \cite{shi2021robust} on networked NNI systems, \cite{Angeli2006} on counter-clockwise (CCW) systems which is closely related to NNI systems, and our conference paper \cite{ghallab2018extending} on single-input-single-output (SISO) NNI systems. However, this paper differs from these works in significant ways: (a) while the paper \cite{shi2021robust} builds on the results of \cite{ghallab2018extending}, this paper establishes the stability robustness of the feedback interconnections of NNI systems under reduced set of assumptions; (b) we generalize here the results of \cite{ghallab2018extending} to multi-inputs-multi-outputs (MIMO) nonlinear negative imaginary systems and, in addition, we investigate the case of NI systems with free motion in the nonlinear domain;
(c) \cite{Angeli2006} investigated the boundedness of the solution of CCW feedback systems whereas we establishes the robust asymptotic stability of NNI feedback systems.

Transforming from the classic frequency domain to time domain for analysis of NNI systems, raises the need for time-domain interpretation of the NNI property. Motivated by the time-domain characterization of the negative imaginary property for LTI systems, we introduce a formal definition for the nonlinear negative imaginary property of a general nonlinear system. The NNI definition will be given in terms of a dissipation inequality with an appropriate work rate. This formulation extends and complements the existing classical dissipativity and passivity framework to allow for more general supply rates which involves deivative of the inputs and outputs. Flexible structures with colocated force actuators and position sensors are dissipative according to this new formulation.

By defining the NNI property in a time-domain dissipativity framework, we seek to extend the above LTI NI stability result to a nonlinear setting. In particular, we generalize the above mentioned NI stability result to establish the stability robustness of a positive feedback interconnection of NNI systems. In order to handle these NI nonlinear systems, a generalization of the DC-loop gain condition \eqref{dcloop} to the nonlinear setting will be developed. Then, under mild technical assumptions in the steady-state, a Lyapunov-based approach and an invariance principle will be employed to guarantee the robust stability of the feedback system. The nonlinear NI stability result will be shown to reduce to the case of a feedback interconnection of MIMO LTI negative imaginary systems, where the plant may have poles on the imaginary axis except at the origin. The applicability of the NNI stability result will be illustrated through an example of a nonlinear mass-spring-damper system with force as input and displacement of the mass as output.

Next, the nonlinear negative imaginary systems theory will be further extended to the case of free motion. It will be shown that, under suitable assumptions, a cascade connection of nonlinear system, affine in the input, and single integrator will lead to a nonlinear negative imaginary system. Robust stability result will be established for a positive feedback system of nonlinear NI systems where the plant is nonlinear NI with integrator and the controller satisfies a astrict notion of the negative imaginary property known as weakly strictly nonlinear NI system property.

The structure of the paper is organized as follows. \textbf{Section 2} reviews some of the related  definitions and results in the NI literature. Then, a characterization of the negative imaginary notion in the nonlinear domain will be given in \textbf{Section 3}. In\textbf{ Section 4}, The formulation of the problem of The main result will be given and the stability robustness of a positive feedback system will be established. \textbf{Section 2} is devoted for further extension of the nonlinear NI framework to the case of free motion and related stability result.




\section{Preliminaries}
For the sake of completeness, we review in this section some of the main related definitions and results in the NI literature for LTI systems, which will lay the foundation of establishment of the NNI systems-theoretic framework.

\subsection{Background for NI LTI Systems}

We consider here the following LTI system:
\begin{align}
\label{eq:xdot1}
&\dot{x}(t) = Ax(t)+Bu(t), \\
\label{eq:y1} &y(t) = Cx(t)+Du(t),
\end{align}
where $A \in \mathbb{R}^{n \times n}, B \in \mathbb{R}^{n \times m}$,  $C\in \mathbb{R}^{m \times n}$, and $D\in \mathbb{R}^{m \times m}$. The system \eqref{eq:xdot1} and \eqref{eq:y1} has the $m\times m$ real-rational proper  transfer function $G(s):=C(sI-A)^{-1}B+D$. The transfer function matrix $G(s)$ is said to be strictly proper if $G(\infty)=D=0$. We will use the notation $\left[\begin{array}{l|l}A & B \\ \hline C & D\end{array}\right]$ to denote the state-space realization of the LTI system \eqref{eq:xdot1}-\eqref{eq:y1}.

The frequency domain characterization of the NI property to LTI system \eqref{eq:xdot1}-\eqref{eq:y1} is given in the following definition.
\begin{definition} \cite{lanzon2008}. 
A square transfer function matrix $G(s)$ is called negative imaginary if the following conditions are satisfied:
\begin{enumerate}
\item[1.] $G(s )$ has no pole at the origin and in $\Re[s]>0$;
\item[2.] For all $\omega >0$ such that $j\omega$ is not a pole of $G(s )$,  $j\omega\left( G(j\omega )-G(j\omega )^{T }\right) \geq 0$;
\item[3.] If $j\omega_{0}$, $\omega_0\in(0,\infty)$ is a pole of $G(j\omega )$, it is at most a simple pole and the residue matrix $K_{0}= \lim_{ s\rightarrow j\omega_{0}}(s-j\omega_{0})sG(s)$ is positive semidefinite Hermitian.
\end{enumerate}
\end{definition}

\noindent The LTI system \eqref{eq:xdot1} and \eqref{eq:y1} is said to be negative imaginary if its transfer function is NI. A stronger version of the negative imaginary property is given in the following definition.


\begin{definition} \cite{lanzon2008}. A square real-rational transfer function matrix $G(s)$ is strictly negative imaginary if:
\begin{enumerate}
\item [1)] $G(s)$ has no poles in $\Re[s]\geq0$;
\item [2)] $j[G(j\omega)-G^{T}(j\omega)]>0$ for $\omega\in(0,\infty)$.
\end{enumerate}
\end{definition}

The next lemma is known as \emph{the negative imaginary lemma} which provides a state-space characterization of LTI NI systems in terms of a linear matrix inequality (LMI).

\begin{lemma}\cite{Mabrok2011}. \label{lem1}
Let $(A,B,C,D)$ be a minimal state-space realization of the transfer function matrix $G(s)$. Then, $G(s)$ is negative imaginary if and only if $det(A)\neq0$, $D=D^T$ and there exist matrices $P=P^T>0$, $W \in \mathbb{R}^{m \times m}$, and $L \in \mathbb{R}^{m \times n}$ such that the following LMI is satisfied:
\begin{equation}\label{rr}
\begin{aligned}
\left[
  \begin{array}{cc}
    PA+A^TP & PB-A^TC^T \\
    B^TP-CA & -(CB+B^TC^T)
  \end{array}
\right]
 &=\left[
  \begin{array}{cc}
    -L^TL & -L^TW \\
    -W^TL & -W^TW
  \end{array}
\right]\\[6pt] &\leq 0.
\end{aligned}
\end{equation}
\end{lemma}

\begin{remark}
 The linear matrix equality \eqref{rr} can be simplified to the following set of equations (see \cite{lanzon2008}),
 $$
AP+PA^{T}\leq 0, \quad \text{and } \quad B+APC^{T}=0.
$$
\end{remark}
The following lemma provides a state-space characterization of strictly negative imaginary LTI systems.
\begin{lemma}(\cite{xiong21010jor}). \label{SNI}
Let $(A,B,C,D)$ be a minimal state-space realization of the transfer function matrix $G(s)$. Then $G(s)$ is strictly negative imaginary if and only if:
\begin{enumerate}
\item [1)] $\det(A)\neq 0$, $D=D^{T}$;
\item [2)] there exists a matrix $P=P^{T}>0, P \in \mathbb{R}^{n \times n}$, such that
$$
AP^{-1}+P^{-1}A^{T}\leq 0, \quad \text{and } \quad B+AP^{-1}C^{T}=0;
$$
\item [3)]the transfer function matrix $M(s)\backsim
\begin{bmatrix}
\begin{array}{c|c}
A & B \\ \hline LPA^{-1} & 0
\end{array}
\end{bmatrix}$
 has full column rank at $s=jw$ for any $\omega\in(0,\infty)$ where $L^{T}L=-AP^{-1}-P^{-1}A^{T}$. That is, rank $M(j\omega)=m$ for all $\omega\in(0,\infty)$.
\end{enumerate}
\end{lemma}


The next lemma provides a time-domain condition equivalent to the frequency-domain definition of the NI property for an LTI system of the form \eqref{eq:xdot1} and \eqref{eq:y1}.

\begin{lemma}\label{time_def} Suppose that the system \eqref{eq:xdot1}-\eqref{eq:y1} (with $D=0$) is controllable and observable.  
Then, the transfer function $G(s)$ is negative imaginary if and only if there exists a matrix $P>0$ such that along the trajectories of the system, the function $V(x)=\frac{1}{2}x^TPx$ satisfies
\begin{equation}\label{lyap}
    \dot{V}(x(t))\leq \dot{y}(t)u(t), \quad \forall \ t\geq0.
\end{equation}
\end{lemma}
\begin{pf}
Differentiating the function $V$ with respect to the time $t$, one has
$$
    \dot{V}(x(t))= \frac{1}{2}x^T(PA+A^TP)x+x^TPBu.
$$
Substituting into the dissipation inequality \eqref{lyap}, we get
\begin{align*}
  \frac{1}{2}x^T(PA+A^TP)x&+x^TPBu\leq  \\
   &u^TCAx+\frac{1}{2}u^T(CB+B^TC^T)u,
\end{align*}
for all $x$ and $u$. In matrix form, the above inequality is equivalent to
$$
\frac{1}{2} \left[\begin{matrix}
x^{T} \hspace{1.5mm} u^{T}
\end{matrix}\right] \left[
  \begin{matrix}
    PA+A^TP & PB-A^TC^T \\
    B^TP-CA & -(CB+B^TC^T) \\
  \end{matrix}
\right] \begin{bmatrix}
x \\
u \\
\end{bmatrix}\leq 0,
$$
for all $x$ and $u$. By Lemma \ref{lem1}, the proof follows.
\begin{flushright}
$\blacksquare$
\end{flushright}
\end{pf}

\subsection{Relation between PR and NI Systems}

An important class of LTI systems of related interest is the class of positive real (PR) systems. The following two lemmas highlight this relation between NI and PR LTI systems in the frequency domain.

\begin{lemma}(\cite{xiong21010jor}). Let $(A,B,C,D)$ be a minimal state-space realization of a transfer function $G(s)$ where $A \in \mathbb{R}^{n \times n},B \in \mathbb{R}^{n \times m},C\in \mathbb{R}^{m \times n},D \in \mathbb{R}^{m \times m}$ and let $\tilde{G}(s):=G(s)-D$. Then $G(s)$ is NI if and only if $F(s):=s\tilde{G}(s)$ is positive real.
\end{lemma}
\begin{lemma}(\cite{buscarino2016forward}).\label{pr_ni} Given a square proper positive real transfer function matrix $G(s)$, then $R(s):=\frac{G(s)}{s}$ is negative imaginary.
\end{lemma}

\subsection{Internal Stability of LTI NI Feedback System}
The stability robustness of positive feedback interconnections of LTI NI systems has been established in \cite{lanzon2008, petersen2010}.
The following lemma will be useful in establishing the Lyapunov-based stability result.
\begin{lemma}(\cite{ghallab2017lyapunov}).\label{lemma1}
Given negative imaginary $G(s)$ with state space realization $(A_1, B_1, C_1, D_1)$ and strictly negative imaginary $H(s)$ with state space realization $(A_2, B_2, C_2, D_2)$. Assume $G(\infty)H(\infty)=0$ and $H(\infty)\geq0$. Also, let $P_1>0$ and $P_2>0$ be corresponding matrices satisfying the LMI \eqref{rr}. Then, the matrix
$$
\left[ \begin{array}{cc}
{P}_{1}-C_{1}^{T}D_{2}C_{1} & -C_{1}^{T}C_{2} \\
{-C}_{2}^{T}C_{1} & {P}_{2}-C_{2}^{T}D_{1}C_{2} \end{array}
\right]
$$
is positive-definite if and only if $\lambda_{max}(G(0)H(0))<1$.
\end{lemma}

One of the main stability results for the positive feedback interconnection of negative imaginary LTI systems is given in the next lemma.
Consider a minimal state-space representation for the NI transfer function $G(s)$,
\begin{align}
\label{eq:xdot11}
&\dot{x}_{1}(t) = A_{1} x_{1}(t)+B_{1} u_{1}(t), \\
\label{eq:y11} &y_{1}(t) = C_{1} x_{1}(t)+D_{1} u_{1}(t),
\end{align}%
where $A_1 \in \mathbb{R}^{n \times n}, B_1 \in \mathbb{R}^{n \times m},C_1
\in \mathbb{R}^{m \times n}, D_1 \in \mathbb{R}^{m \times m}$.
\noindent Also, we consider a minimal state-space representation for the SNI transfer function $H(s)$,
\begin{align}
\label{eq:xdot2}
&\dot{x}_{2}(t) = A_{2} x_{2}(t)+B_{2} u_{2}(t), \\
\label{eq:y2} &y_{2}(t) = C_{2} x_{2}(t)+D_{2} u_{2}(t),
\end{align}%
where $A_2 \in \mathbb{R}^{n \times n}, B_2 \in \mathbb{R}^{n \times m}, C_2
\in \mathbb{R}^{m \times n},D_2 \in \mathbb{R}^{m \times m}$.

\begin{theorem}(\cite{ghallab2017lyapunov}). \label{NI-stability}
Assume that $G(s)$ is negative imaginary with minimal realization \eqref{eq:xdot11} and \eqref{eq:y11} and $H(s)$ is strictly negative imaginary with minimal realization \eqref{eq:xdot2} and \eqref{eq:y2} such that $G(\infty)H(\infty)=0$ and $H(\infty)\geq0$. Also, assume that $\lambda_{max}(G(0)H(0))<1$. Then, the positive feedback interconnection of $G(s)$ and $H(s)$ as in Figure \ref{sys_31} is internally stable.
\end{theorem}
The proof of the above theorem has been established in \cite{ghallab2017lyapunov} where a Lyapunov based technique has been employed to establish the internal stability of the positive feedback interconnection of $G(s)$ and $H(s)$. With $V_1(x_1)=\frac{1}{2}x_{1}^{T}P_1x_1$ and $V_2(x_2)=\frac{1}{2}x_{2}^{T}P_2x_2$, the function
\begin{align*}
 V(x_{1},x_{2}):&=V_1(x_1)+V_2(x_2)-2y_{1}^{T}y_{2}\\
                &=\frac{1}{2}x_{1}^{T}P_1x_1+\frac{1}{2}x_{2}^{T}P_2x_2-2y_{1}^{T}y_{2}
\end{align*}
%
is considered as Lyapunov function candidate in terms of the trajectories of the feedback subsystems. Rewriting this function we obtain
\begin{align*}
  V&(x_1,x_2)\\
  &=\frac{1}{2}x_1^T(P_1-C_1^TD_2C_1)x_1+\frac{1}{2}x_2^T(P_2-C_2^TD_1C_2)x_2\\
            &\hspace{1.5cm}-x_1^TC_1^TC_2x_2\\[8pt]
            &= \frac{1}{2} \left[\begin{array}{r}
                x_{1}^{T} \hspace{1.5mm} x_{2}^{T}
               \end{array}\right]
               \left[ \begin{array}{cc}
               {P}_{1}-C_{1}^{T}D_{2}C_{1} & -C_{1}^{T}C_{2} \\
               {-C}_{2}^{T}C_{1} & {P}_{2}-C_{2}^{T}D_{1}C_{2} \end{array}
                \right]
               \left[ \begin{array}{c}
                  x_{1} \\
                  x_{2} \\
               \end{array}\right]
\end{align*}
By Lemma \ref{lemma1}, it follows the function $V(x_1,x_2)$ is positive-definite and thus is a valid candidate Lyapunov function. The time-derivative of $V(x_1,x_2)$ has been shown to be negative semidefinite in \cite{ghallab2017lyapunov}, and hence the feedback system of $G(s)$ and $H(s)$ is at least Lyapunov stable. Then, using the assumptions of Theorem \ref{NI-stability}, all the poles of the closed-loop system matrix has been shown in \cite{ghallab2017lyapunov} to lie in the open left-half of the complex plane, which proves the internal stability of the positive feedback interconnection of $G(s)$ and $H(s)$.

\section{Nonlinear Negative Imaginary Systems}

In this section, we will seek to establish a natural generalization of the negative imaginary property of an LTI to general nonlinear systems. In particular, a nonlinear system which is passive from the input to the derivative of the output will be termed as a \emph{nonlinear negative imaginary (NNI) system }. Also, later in this section, we seek to establish a stability robustness result for  NNI feedback systems analogous to the passivity theorem. First of all, we provide examples that motivate the need for a system-theoretic framework that extends classical dissipativity\slash passivity theory to allow for systems with supply rates involve the derivative of the input and outputs of the system.

\subsection{Examples of NNI Systems}



\begin{example}[Mass-Spring-Damper System]\label{msd_exa}
\leavevmode

\normalfont
For mechanical systems, any mass-spring-damper system with non-negative damping coefficients is passive with respect to force inputs and corresponding velocity outputs; i.e., nonlinear negative imaginary from the force input and displacement output. We claim that the system has nonlinear negative imaginary I-O dynamics from $u$ to the output $x$. To prove this claim, we consider the following equation for the dynamics of a nonlinear mass-spring-damper system:
\begin{equation}
\ddot{x}+\beta(x, \dot{x}) \dot{x}+k(x)=u, \qquad y=x
\end{equation}
where $m$ is the mass of the block attached to the spring, $\beta(x,\dot{x})\geq 0$ is a dynamic state-dependent friction coefficient and the term $k(x)$ is the nonlinear spring stiffness such that $k(0)=0$. In fact, considering the total energy of the system as a storage function, we get
\begin{equation}
V(x(t),\dot{x}(t))=\frac{1}{2}m\dot{x}^2+\int_{0}^{x}k(\xi)d\xi.
\end{equation}
Taking the time derivative of $V$ implies
\begin{align*}
  \dot{V}(x(t),\dot{x}(t))&=m\dot{x}\ddot{x}+k(x)\dot{x}=\dot{x}u-\beta(x,\dot{x})\dot{x}^2\\
                          &\leq \dot{x} u.
\end{align*}
This completes the proof.
\begin{flushright}
$\blacksquare$
\end{flushright}
\end{example}

\begin{example}[Hamiltonian Systems]
\leavevmode

\normalfont
A Hamiltonian system which is affine in the input, has negative-imaginary I-O dynamics from the input force $F$ to the corresponding coordinate $q$. To establish this claim, consider  a Hamiltonian system which satisfies the following equations:
\begin{align*}
\dot{q}&=\frac{\partial H}{\partial p}(q, p); \\
\dot{p}&=-\frac{\partial H}{\partial q}(q, p)+F; \\
      y&=q.
\end{align*}
Choosing a storage function $V$ as the Hamiltonian $H$ of the system and taking the time derivative of $V$ yields:
\begin{equation}
  \dot{V}=\dot{H}=\frac{\partial^T H}{\partial q} \dot{q}+\frac{\partial^T H}{\partial p} \dot{p}=\frac{\partial H}{\partial^T p} u=\dot{y}^{T} F
\end{equation}
which completes the proof.
\begin{flushright}
$\blacksquare$
\end{flushright}
\end{example}
This NNI property can also be established for general Euler-Lagrangian systems as in the following example.
\begin{example}[Euler-Lagrangian Systems]
\leavevmode

\normalfont
Many standard mechanical and physical systems are modeled by Euler-Lagrangian dynamics as follows:
\begin{equation}\label{lagrange}
{M}({q}) \ddot{{q}}+{C}({q}, \dot{{q}})\dot{{q}}+{{G(q)}}={F},
\end{equation}
where ${M(q)}$ is a positive-definite symmetric inertia matrix. The term ${C({q}, \dot{q})}$ is the Coriolis and centrifugal matrix where ${C}({{q}}, \dot{{q}})=\frac{\mathrm{d}}{\mathrm{d} t}{M(q)}-\frac{1}{2} \frac{\partial}{\partial {q}}\left(\dot{{q}}^{T} {M} \right)$, ${G(q)}$ is the gravitational vector, and ${F}$ is the generalized input of the  system. The matrix ${C}({{q}}, \dot{{q}})$ is normally defined using the Christoffel symbols; in which ${\dot{M}} - 2{C}({{q}}, {\dot{q}})$ is skew-symmetric \cite{ortega1989adaptive, lewis2003robot}.

To show that the system \eqref{lagrange}is NNI, we consider the following storage function given by
\begin{equation}
    V({q}, \dot{{q}})=\frac{1}{2}\dot{{q}}^{T} {M}({q}) \dot{{q}}+{P(q)}.
\end{equation}
Here, ${P(q)}$ is the potential vector, where ${G(q)}=\frac{\partial {P(q)}}{\partial {q}}$, and it is assumed to have an absolute minimum at ${q}={0}$. It can be easily shown that $V$ is positive-definite since ${M(q)}$ is a positive-definite matrix, and ${P(q)}$ is a nonnegative scalar quantity. Then, taking the first time derivative of the function $V$ we obtain
\begin{equation*}
\begin{aligned}
\frac{d V}{d t}({q}, {\dot{q}}) &={\dot{q}}^{T} {M}({q}) {\ddot{q}}+\frac{1}{2} {\dot{q}}^{T} {\dot{M}}({q}) {\dot{q}}+{G({q})}{\dot{q}} \\
&={\dot{q}}^{T}\left[-{C}({q}, {\dot{q}}) {\dot{q}}-{G(q)}+{F}\right]\\ & \hspace{2cm} +\frac{1}{2} {\dot{q}}^{T} {\dot{M}}({q}) {\dot{q}}+{G(q)}{\dot{q}} \\
&={\dot{q}}^{T} {F}+\frac{1}{2} {\dot{q}}^{T}\left[{\dot{M}}({q})-2 {C}({q}, {\dot{q}})\right] {\dot{q}}= {\dot{q}}^{T} {F},
\end{aligned}
\end{equation*}
which shows that the above Euler-Lagrangian system \eqref{lagrange} is passive system from ${F}\mapsto \dot{q}$; i.e., nonlinear negative imaginary from ${F}\mapsto {q}$. It is worth noting there is no notion of passivity for the map ${F}\mapsto{q}$. In other words, for systems represented by \eqref{lagrange}, we will be able to prove passivity from the control input $F$ to the velocity of the generalized coordinates $\dot{q}$ but not to the position of the generalized coordinates $q$.
\begin{flushright}
$\blacksquare$
\end{flushright}
\end{example}

\subsection{Characterization of NNI Systems}

Motivated by the time-domain characterization of the negative imaginary property for LTI systems (Lemma \ref{time_def}), we introduce here a formal definition of the nonlinear negative imaginary property for a general nonlinear system in terms of a dissipativity inequality with an appropriate work rate involving the derivative of the output of the system.

Consider the following general MIMO nonlinear system of the form
\begin{align}
\label{x}\dot{x}&= f(x,u),\\
\label{y} y&= h(x)
\end{align}
where $f:\mathbb{R}^n\times \mathbb{R}^m\rightarrow\mathbb{R}^n$ is Lipschitz continuous function such that $f(0,0)=0$ and $h:\mathbb{R}^n\rightarrow\mathbb{R}^m$ is continuously differentiable function such that $h(0)=0$. We have the following definition.

\begin{definition}\label{nonlinear_NI}\normalfont The system \eqref{x}-\eqref{y} is said to be \textbf{\emph{nonlinear negative imaginary}} if there exists a positive-definite  function $V:\mathbb{R}^n\rightarrow\mathbb{R}$ (called a storage function) of a class $C^1$ such that the following dissipative inequality
\begin{equation}\label{lyapD}
    \dot{V}(x(t))\leq \dot{y}^T(t)u(t),
\end{equation}
holds for all $t \geq 0$.
\end{definition}
The dissipative inequality can be given alternatively by integrating the infinitesimal form \eqref{lyapD} from $0$ to $t$ to get its integral equivalent form; i.e,
\begin{equation}\label{lyapI}
  V(x(t))\leq V(x(t_0))+\int_{0}^{t}\dot{y}^T(s)u(s)ds
\end{equation}
for all $t \geq 0$.


The following lemma gives necessary and sufficient conditions for a nonlinear system of the form \eqref{x}, \eqref{y}, affine in the input, to be nonlinear negative imaginary with a $C^1$ storage function, and is deemed analogous to the nonlinear KYP lemma for passive systems \cite{isidori1999asymptotic}.

\begin{lemma}
Consider a system of the form \eqref{x}-\eqref{y} and affine in the control, that is described by equations of the form
\begin{equation}\label{kyp_affine}
\begin{aligned}
  \dot{x}(t)=& f(x(t))+g(x(t))u(t)\\
  y(t)=& h(x(t)),
\end{aligned}
\end{equation}
where $f, g, h$ are smooth functions of $x$ with $f(0)=0$ and $h(0)=0$. Then the system \eqref{kyp_affine} is NNI~ if and only if there exist functions $V(x)\geq0$, $L(x)$, and $W(x)$ such that
\begin{equation}\label{kyp_affine1}
\begin{aligned}
  \nabla V^{T}f(x) &= -L^{T}(x)L(x) \\
  \frac{1}{2}g^{T}(x)\nabla V(x)&=\frac{1}{2}\nabla h(x)f(x) -W^{T}(x)L(x)\\
  W^{T}(x)W(x)&= \frac{1}{2}(\nabla h(x)g(x))^{T}+\nabla h(x)g(x)
\end{aligned}
\end{equation}
\end{lemma}
\begin{pf}\normalfont
(\textbf{Necessity}) From the dissipativity inequality \eqref{lyapD}, for all $x$ and $u$, we get
\begin{equation*}
\nabla V^T(x) f(x) +\nabla V^T g(x)u \leq (f^T(x)+u^Tg^T(x))\nabla h^T(x) u.
\end{equation*}
Rearranging the terms of the above inequality, we obtain
\begin{align*}
    f^T(x)\nabla h^T(x) u &+u^Tg^T(x)\nabla h^T(x) u \\ &- \nabla V^T(x) f(x) +\nabla V^T(x) g(x)u\geq0,
\end{align*}
that is; the right-hand side of the above inequality can be factored, for all $x$ and $u$, as follows
\begin{align*}
f^T&(x)\nabla h^T(x) u +u^Tg^T(x)\nabla h^T(x) u - \nabla V^T(x) f(x)\\
                &+\nabla V^T(x) g(x)u=(L(x)+W(x)u)^T(L(x)+W(x)u)\\
                &=L^T(x)L(x)+2L^T(x)W(x)u+u^TW^T(x)W(x)u
\end{align*}
for some $L(x)$ and $W(x)$. By equating coefficients of similar powers of $u$, we get the equations \eqref{kyp_affine1}.

\textbf{(Sufficiency)} For all $x$ and $u$, we have
\begin{align*}
  \dot{y}^T(t)u(t)&= (\nabla h(x).(f(x)+g(x)u(t)))^T u(t) \\
                  &= f^T(x)\nabla h^T(x) u +u^Tg^T(x)\nabla h^T(x) u\\
                  &= \nabla V^T(x) g(x)u + W^T(x)L(x) u\\ & +u^TW^T(x)W(x)u +\nabla V^T(x) f(x)+L^T(x)L(x)\\
                  &= \nabla V^T(x) (f(x)+g(x)u)\\& \hspace{1.5cm}+(L(x)+W(x)u)^T(L(x)+W(x)u)\\
                  &\leq \nabla V^T(x)\dot{x}= \dot{V}(x(t)),
\end{align*}
which concludes the proof.
\begin{flushright}
$\blacksquare$
\end{flushright}
\end{pf}
%

\subsection{Interconnections of NNI Systems}

Similar to passive systems, the class of systems with NNI dynamics is closed under parallel and feedback interconnections. This fact can be readily verified if we consider the storage function $V$ of the interconnected system as the sum of the storage functions of the subsystems storage functions $V_1$ and $V_2$; that is
\begin{equation*}
 V(x) = V_1(x_1) + V_2(x_2), \quad x=(x_1^T,x_2^T)^T
\end{equation*}
For example, the parallel interconnection of NNI systems is characterized by $u=u_1 = u_2 $, $y = y_1 + y_2$. Then the time derivative of the storage function $V$ is evaluated as
\begin{equation*}
\dot{V}=\dot{V}_1+\dot{V}_2\leq \dot{y}_1^Tu_1+\dot{y}_2^Tu_2=(\dot{y}_1+\dot{y}_2)^T u=\dot{y}^T u
\end{equation*}
Similarly we can show that positive feedback interconnection of two nonlinear NI systems is again a nonlinear NI system.  The proof is straightforward recalling the equation $u_1 = y_2$, $u_2 = y_1 $ and proceeding as described above.
\subsection{Stronger Notions of NNI}
We conclude this section by introducing some stronger notions of the NNI property in an analogous way \cite{isidori1999asymptotic} for passive systems. These notions will be employed in the next section for the stability analysis of NNI feedback systems in a similar fashion to that of passive systems. 



\begin{definition}\normalfont
The system \eqref{x}-\eqref{y} is said to be \textbf{\emph{marginally strictly nonlinear negative imaginary (MS-NNI)}} if the dissipativity inequality \eqref{lyapD} is satisfied, and in addition, if $u$, $x$ are such that
\begin{equation}\label{06}
   \dot{V}(x)=\dot{y}^T(t)u(t)
\end{equation}
for all $t>0$, then $\lim_{t\rightarrow\infty}u(t)=0$.
\end{definition}
\begin{definition}\label{wsnni}\normalfont
The system \eqref{x}-\eqref{y} is said to be \textbf{\emph{weakly strictly nonlinear negative imaginary (WS-NNI)}} if it is MS-NNI and globally asymptotically stable when $u\equiv0$.
\end{definition}
\begin{remark}
The above definition of weakly strictly nonlinear negative imaginary systems is deemed a nonlinear analog to the SNI property of LTI systems. In Appendix A, we show that for an LTI system of the form \eqref{eq:xdot1} and \eqref{eq:y1}, the WS-NNI property reduces to the SNI property.
\end{remark}

\section{Stability of Interconnected Nonlinear Negative Imaginary Systems}
%
We aim in this section to obtain a result on the asymptotic stability of positive feedback interconnections of nonlinear negative imaginary systems. To attain this objective, we seek to establish a nonlinear generalization of the negative imaginary lemma. Then, Lyapunov stability theory and LaSalle's invariance principle will be employed to establish asymptotic stability.

To this end, we consider the following general MIMO nonlinear systems:
\begin{align*}\label{sys}
  H_1{:} \quad {\dot{{x}}}_1 & = {f}_1({x}_1, {u}_1),\\
                                {y}_1 & = {h}_1({x}_1)  \\
                        \intertext{and}
 H_2{:} \quad  {\dot{{x}}}_2 & =  {f}_2({x}_2, {u}_2),\\
                                {y}_2 & = {h}_2({x}_2)\\
\end{align*}
where ${h}_i:\mathbb{R}^n\rightarrow\mathbb{R}^n$ is a ${C}^1$ function with ${h}_i({0})={0}$, ${f}_i:\mathbb{R}^n\times \mathbb{R}^n\rightarrow\mathbb{R}^n$ is continuous and locally Lipschitz in ${x}_i$ for bounded ${u}_i$, and where ${f}_i({0},{0})={0}$. We first construct a candidate Lyapunov function for the interconnection of systems $H_1$ and $H_2$. This interconnected system determines the stability properties of the closed-loop system, see \cite{ghallab2018extending}.


\subsection{Open-Loop System Result}

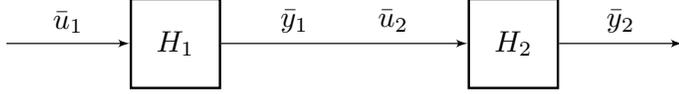
\begin{figure}[t]
\centering
\tikzstyle{block} = [draw, thick, rectangle,
    minimum height=2em, minimum width=4em]
\tikzstyle{sum} = [draw, circle,inner sep=0pt,minimum size=1pt, node distance=1cm]
\tikzstyle{input} = [coordinate]
\tikzstyle{output} = [coordinate]
\tikzstyle{pinstyle} = [pin edge={to-,thin,black}]
\tikzstyle{int}=[draw, thick, rectangle, minimum height = 3em,
    minimum width = 3em]
\resizebox{9cm}{!}{%
\begin{tikzpicture}[node distance=2.5cm,auto,>=latex']
\tikzstyle{block} = [draw, thick, rectangle,
    minimum height=3em, minimum width=6em]
    \node [int] (a) {$H_1$};
    \node (b) [left of=a,node distance=2cm, coordinate] {a};
    \node [int] (c) [right of=a,  node distance=4cm] {$H_2$};
    \node [coordinate] (end) [right of=c, node distance=2cm]{};
    \path[->] (b) edge node {$\bar{u}_1$} (a);
    \path[->] (a) edge node {$\bar{y}_1  \quad \quad  \ \bar{u}_2$} (c);
    \draw[->] (c) edge node {$\bar{y}_2$} (end) ;
\end{tikzpicture}}
\caption{Open-loop interconnection of systems $H_1$ and $H_2$ in steady-state.}\label{sys13}
\end{figure}

Inspired by the result of Lemma, we shall consider the open-loop interconnection of systems $H_1$ and $H_2$. We seek to construct a candidate Lyapunov function using a set of technical assumptions on the systems $H_1$ and $H_2$ and on the open-loop interconnection of the systems $H_1$ and $H_2$ as shown in Figure \ref{sys13}. These assumptions are a nonlinear generalization of the assumptions of Lemma.

We have the following assumptions for the systems of $H_1$ and $H_2$ and their the open-loop interconnection.
\begin{assumption}\label{I}
    For any constant $\bar{{u}}_1$, there exists a unique solution $(\bar{{x}}_1, \bar{{y}}_1)$ to the equations
\begin{align*}
  {0}  &= {f}_1(\bar{{x}}_1,\bar{{u}}_1),\\
  \bar{{y}}_1&= {h}_1(\bar{{x}}_1)
\end{align*}
such that $\bar{{u}}_1\neq {0}$ implies $\bar{{x}}_1\neq {0}$ and the mapping $\bar{{u}}_1 \mapsto \bar{{x}}_1$ is continuous.
\end{assumption}
\vspace{.25cm}
\begin{assumption}\label{II}
    For any constant $\bar{{u}}_2$, there exists a unique solution $(\bar{{x}}_2, \bar{{y}}_2)$ to the equations
\begin{align*}
 {0}  &= {f}_2(\bar{{x}}_2,\bar{{u}}_2),\\  \bar{{y}}_2&= {h}_2(\bar{{x}}_2)
\end{align*}
such that $\bar{{u}}_2\neq {0}$ implies $\bar{{x}}_2\neq {0}$.
\end{assumption}
\vspace{.25cm}
\begin{assumption}\label{III}
   ${h}_1^T(\bar{{x}}_{1}){h}_2(\bar{{x}}_{2})\geq 0$, for any constant $\bar{{u}}_1$ where $\bar{{u}}_2=\bar{y}_1$.
    \vspace{.25cm}
\end{assumption}
\begin{assumption}\label{IV}
    For any constant $\bar{{u}}_1$, let $(\bar{{x}}_1, \bar{{y}}_1)$ be defined as in Assumption \ref{I}. Also, let $(\bar{{x}}_2, \bar{{y}}_2)$ be defined as in Assumption \ref{II} where $\bar{{u}}_2=\bar{{y}}_1$. Then there exits a constant $0<\gamma<1$ such that for any $\bar{{u}}_1$, the following sector bound condition:
  \begin{equation}\label{sec}
  \bar{{y}}_2^T\bar{{y}}_2\leq \gamma^2 \bar{{u}}_1^T\bar{{u}}_1.
\end{equation}
\end{assumption}
\vspace{1em}
Next we state a lemma which establishes a candidate Lyapunov function for the positive feedback interconnection of the systems $H_1$ and $H_2$. This lemma is a nonlinear analog of Lemma~\ref{lemma1}.
\begin{lemma}\label{lemma12}
Assume that the system $H_1$ is NNI with corresponding storage function $V_1(x_1)$ and $H_2$ is WS-NNI with corresponding storage function $V_2(x_2)$. Consider the function defined by
\begin{equation}\label{lyap_fn}
W({x}_1,{x}_2):=V_1({x}_1)+V_2({x}_2)-{h}_1^T({x}_1){h}_2({x}_2).
\end{equation}
Suppose that the Assumptions \ref{I}-\ref{IV} are satisfied. Then the function $W$ is positive-definite.
\end{lemma}
\begin{pf} \normalfont
We have $W({0},{0})=0$. Next, we show that the function $W({x}_1,{x}_2)>0$ for all nonzero $(x_1,x_2)$. The proof will be established using an iterative technique and the fact that for any constant input ${u}_i(t)\equiv {\bar{u}}_i$ each subsystem satisfies
\begin{equation}\label{modlyapd}
  V_i({x}_i(t))-{h}_i^T(x_i(t)){\bar{u}}_i\geq V_i({\xi}_i)-{h}_i^T({\xi}_i){\bar{u}}_i,
\end{equation}
for all ${x}_i, \ {\xi}_i \in\mathbb{R}^n$. This fact can be readily seen by integrating the following dissipativity inequality
\begin{equation}\label{non_NI}
  \dot{V_i}({x}_i(t))\leq \dot{{y}}_i^T(t){{\bar{u}}}_i.
\end{equation}
from $0$ to $t$. Now fix any $({x}_1, {x}_2)\neq{0}$, and consider the iteration $\big\{\bar{{u}}_1^{[i]}, \bar{{x}}_1^{[i]}, \bar{{y}}_1^{[i]}, \bar{{y}}_2^{[i]}\big\}$ such that
\begin{align*}
   \bar{{u}}_1^{[1]} &:= {h}_2({x}_2), \ \text{and}\\[4pt]
   \bar{{u}}_1^{[i+1]}&:=\frac{1}{\sqrt{\upgamma}}\bar{{y}}_2^{[i]}, \ i=1,2,3,\ldots
\end{align*}
\noindent We have the following two cases:

\textbf{\textbf{Case $\mathbf{1)}$}} if ${h}_2({x}_2)={0}$, then we have
\begin{equation*}
W({x}_1,{x}_2)=V_1({x}_1)+V_2({x}_2)>0,
\end{equation*}
since $({x}_1,{x}_2)\neq0$ and $V_1(\cdot), V_2(\cdot)$ are positive-definite.

\textbf{Case} $\mathbf{2)}$ if $\bar{{u}}_2={h}_2({x}_2)\neq{0}$, then the condition \eqref{sec} implies
\begin{equation}
\|{\bar{y}}_2\|\leq \upgamma ~ \|{\bar{u}}_1\|,
\end{equation}
and thus from the above iteration we have
\begin{equation*}
\|\bar{{u}}_1^{[i+1]}\|= \ddfrac{1}{\sqrt{\upgamma}} \|\bar{{y}}_1^{[i]}\|\leq \frac{\upgamma}{\sqrt{\upgamma}} \|\bar{{u}}_1^{[i]}\|  , \quad i=1,2,\ldots
\end{equation*}
that is
\begin{equation*}
\|\bar{{u}}_1^{[i+1]}\|\leq \sqrt{\upgamma} ~ \|\bar{{u}}_1^{[i]}\|,
\end{equation*}

which implies $\bar{{u}}_1^i\rightarrow {0} \ (\text{in } \mathbb{R}^n) \ \text{as } i \rightarrow \infty$. Also, according to Assumptions \ref{I} and \ref{II}, for each $i$ the input $ \bar{{u}}_1^{[i]}$ guarantees the corresponding quantities $\bar{{x}}_1^{[i]},\bar{{y}}_1^{[i]},\bar{{u}}_2^{[i]},\bar{{x}}_2^{[i]},\bar{{y}}_2^{[i]}$ and therefore
$$
  \bar{{x}}_1^{[i]}, \bar{{y}}_1^{[i]}, \bar{{u}}_2^{[i]}, \bar{{x}}_2^{[i]}, \bar{{y}}_2^{[i]}\rightarrow {0}, \quad \text{as  \  } i \rightarrow \infty.
$$
\begin{figure}[t]
  \centering
  \includegraphics[scale=0.9]{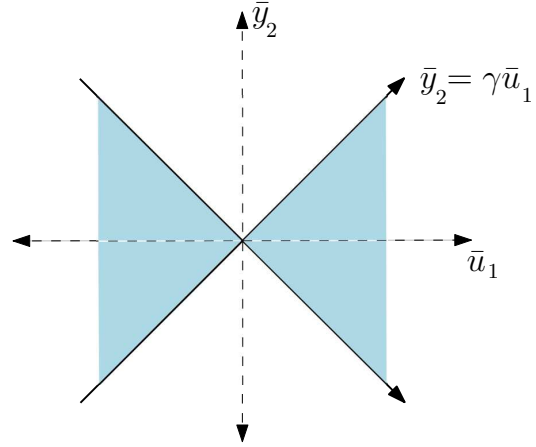}
  \caption{Graphical representation of the sector bound condition \eqref{sec}}\label{sector_condition}
\end{figure}
\noindent Now by exploiting inequality \eqref{modlyapd} we seek to find a lower bound for the function $W(x_1,x_2)$ where $({x}_1, {x}_2)\neq{0}$.

\noindent Let ${\bar{u}}_1=\bar{{u}}_1^{[1]}:=h_2({x}_2)$ and ${\bar{\xi}}_1=\bar{{x}}_1^{[1]}$ in \eqref{modlyapd} for the system $H_1$. Then we obtain
$$
V_1({x}_1)-{h}_1^T({x}_1){h}_2({x}_2)\geq V_1(\bar{{x}}_1^{[1]})-{h}_1^T(\bar{{x}}_1^{[1]}){h}_2({x}_2).
$$
It then follows
$$
  W({x}_1,{x}_2) \geq V_1(\bar{{x}}_1^{[1]})-{h}_1^T(\bar{{x}}_1^{[1]}){h}_2({x}_2)+V_2({x}_2).
$$
Let ${\bar{u}}_2=\bar{u}_2^{[1]}={h}_1(\bar{{x}}_1^{[1]})$ and ${\bar{\xi}}_2=\bar{{x}}_2^{[1]}$ in \eqref{modlyapd} for the system $H_2$. Then we get
$$
V_2({x}_2)-{h}_2^T({x}_2){h}_1(\bar{{x}}_1^{[1]})\geq V_2(\bar{{x}}_2^{[1]})-{h}_2^T(\bar{{x}}_2^{[1]}){h}_1(\bar{{x}}_1^{[1]}).
$$
\noindent This leads to
$$
  W({x}_1,{x}_2) \geq V_1(\bar{{x}}_1^{[1]})+ V_2(\bar{{x}}_2^{[1]})-{h}_1^T(\bar{{x}}_1^{[1]}){h}_2(\bar{{x}}_2^{[1]}).
$$
From Assumption \ref{III}, we consider the following three cases:

$\mathbf{\RomanNumeralCaps{1})}$  ${h}_1(\bar{{x}}_1^{[1]})={0}$ and ${h}_2(\bar{{x}}_2^{[1]})\neq {0}$. 
This implies
\begin{equation*}
W({x}_1,{x}_2)\geq V_1(\bar{{x}}_1^{[1]})+V_2(\bar{{x}}_2^{[1]})\geq V_1(\bar{{x}}_1^{[1]})>0,
\end{equation*}
since we have $\bar{{x}}_1^{[1]}\neq {0}$.

$\mathbf{\RomanNumeralCaps{2})}$ ${h}_1(\bar{{x}}_1^{[1]})\neq {0}$ and ${h}_2(\bar{{x}}_2^{[1]})={0}$. In this case,
$$
W({x}_1,{x}_2)\geq V_1(\bar{{x}}_1^{[1]})+V_2(\bar{{x}}_2^{[1]})\geq V_2(\bar{{x}}_2^{[1]})>0, \text{as\ } \bar{{x}}_2^{[1]}\neq{0}.\\[2pt]
$$

$\mathbf{\RomanNumeralCaps{3})}$  ${h}_1^T(\bar{{x}}_{1}^{[1]}){h}_2(\bar{{x}}_{2}^{[1]})> 0$. In this case,
\begin{align}
  \nonumber W({x}_1,{x}_2)  & \geq V_1(\bar{{x}}_1^{[1]})+ V_2(\bar{{x}}_2^{[1]})-{h}_1^T(\bar{{x}}_1^{[1]}){h}_2(\bar{{x}}_2^{[1]})\\[6pt]
         \nonumber     & = V_1(\bar{{x}}_1^{[1]})+ V_2(\bar{{x}}_2^{[1]})-\ddfrac{1}{\sqrt{\upgamma}}{h}_1^T(\bar{{x}}_1^{[1]}){h}_2(\bar{{x}}_2^{[1]})\\
              & \quad +(\ddfrac{1}{\sqrt{\upgamma}}-1){h}_1^T(\bar{{x}}_1^{[1]}){h}_2(\bar{{x}}_2^{[1]}).
\end{align}

Let ${u}_1=\bar{{u}}_1^{[2]}=\ddfrac{1}{\sqrt{\upgamma}}{h}_2(\bar{{x}}_2^{[1]})$ and ${\xi}_1=\bar{{x}}_1^{[2]}$ in \eqref{modlyapd} for the system $H_1$. Then we get
\begin{align*}
  V_1(\bar{{x}}_1^{[1]})-\ddfrac{1}{\sqrt{\upgamma}}{h}_1^T(\bar{{x}}_1^{[1]})&{h}_2(\bar{{x}}_2^{[1]})\geq \\
   & V_1(\bar{{x}}_1^{[2]})-\ddfrac{1}{\sqrt{\upgamma}} {h}_1^T(\bar{{x}}_1^{[2]}){h}_2(\bar{{x}}_2^{[1]}),
\end{align*}
which in turn leads to
\begin{align*}
 W({x}_1,{x}_2)& \geq V_1(\bar{{x}}_1^{[2]})-\frac{1}{\sqrt{\upgamma}}{h}_1^T(\bar{{x}}_1^{[2]}){h}_2(\bar{{x}}_2^{[1]})+V_2(\bar{{x}}_2^{[1]})\\
 &\quad+\big(\ddfrac{1}{\sqrt{\upgamma}}-1\big){h}_1^T(\bar{{x}}_1^{[1]}){h}_2(\bar{{x}}_2^{[1]}).
\end{align*}
Let ${u}_2=\bar{{u}}_2^{[2]}= \ddfrac{1}{\sqrt{\upgamma}}{h}_1(\bar{x}_1^{[2]})$ and ${\xi}_2=\bar{{x}}_2^{[2]}$ in \eqref{modlyapd} for the system $H_2$. Then we get
\begin{align*}
V_2(\bar{{x}}_2^{[1]})-\ddfrac{1}{\sqrt{\upgamma}}{h}_2^T(\bar{{x}}_2^{[1]}){h}_1(\bar{{x}}_1^{[2]})\geq V&_2(\bar{{x}}_2^{[2]})\\ & - \ddfrac{1}{\sqrt{\upgamma}}{h}_1^T(\bar{{x}}_1^{[2]}){h}_2(\bar{{x}}_2^{[2]}),
\end{align*}
which leads to
\begin{align*}
   W({x}_1,{x}_2) & \geq V_1(\bar{{x}}_1^{[2]})+ V_2(\bar{{x}}_2^{[2]})-\ddfrac{1}{\sqrt{\upgamma}}{h}_1^T(\bar{{x}}_1^{[2]}){h}_2(\bar{{x}}_2^{[2]})\\
              & \quad +\big(\ddfrac{1}{\sqrt{\upgamma}}-1\big){h}_1^T(\bar{x}_1^{[1]}){h}_2(\bar{{x}}_2^{[1]}).
\end{align*}
\noindent Repeating the above process, we obtain
\begin{align*}
 W({x}_1,{x}_2)   & \geq V_1(\bar{{x}}_1^{[i]})+ V_2(\bar{{x}}_2^{[i]})-\ddfrac{1}{\sqrt{\upgamma}}{h}_1^T(\bar{{x}}_1^{[i]}){h}_2(\bar{{x}}_2^{[i]})\\
              & \quad +(\ddfrac{1}{\sqrt{\upgamma}}-1){h}_1^T(\bar{{x}}_1^{[1]}){h}_2(\bar{{x}}_2^{[1]}).
\end{align*}
Letting $i \rightarrow \infty$, we conclude that
\begin{align*}
 W({x}_1,{x}_2)  &  \geq   V_1({0})+ V_2({0})-\ddfrac{1}{\sqrt{\upgamma}}{h}_1^T({0}){h}_2({0})\\
              & \quad +(\ddfrac{1}{\sqrt{\upgamma}}-1){h}_1^T(\bar{{x}}_1^{[1]}){h}_2(\bar{{x}}_2^{[1]})\\
              &  =   0+(\ddfrac{1}{\sqrt{\upgamma}}-1){h}_1^T(\bar{{x}}_1^{[1]}){h}_2(\bar{{x}}_2^{[1]})>0.
\end{align*}
Therefore, the function $W$ is positive-definite for all nonzero ${x}_1, {x}_2$. This completes the proof.
\begin{flushright}
$\blacksquare$
\end{flushright}
\end{pf}

\begin{remark}
We point out here that we used in this paper a stronger version of the sector bound condition than that used in \cite{ghallab2018extending}. This is to ensure the convergence of the iterations used in the proof of the above lemma.
\end{remark}

\begin{remark}
In the LTI case, the aforementioned Assumptions \ref{I}-\ref{IV} reduce to a lightly stronger version of conditions of Lemma \ref{lemma1}. To see this, for the open loop systems as in Figure~\ref{sys13}, consider the SISO LTI case where $H_1$ and $H_2$ are two NI systems represented by
\begin{align*}
 H_i{:} \quad  {\dot{x}}_i&=A_ix_i+B_iu_i\\
                                y_i&=C_ix_i, \qquad i=1,2
\end{align*}
where $H_1(s)$ and $H_2(s)$ are the transfer functions for systems $H_1$ and $H_2$, respectively.  We can see that Assumptions \ref{I}, \ref{II} hold trivially for linear systems. Also, Assumption \ref{III} amounts to the condition $H_2(\infty)>0$. The sector bound condition \eqref{sec} reduces to a stronger version of the DC-gain condition \eqref{dcloop}, which can be seen from the following
\begin{align*}
 \bar{y}_2^2\leq \gamma \bar{u}_1^2 &\Rightarrow \bar{y}_2< \bar{u}_1^2 \quad \text{as \ } \gamma<1 \\
   &\Rightarrow \frac{\bar{y}_2^2}{\bar{u}_2^2}<\frac{\bar{u}_1^2}{\bar{y}_1^2} \quad \text{since $\bar{y}_1=\bar{u}_2$}\\
   & \Rightarrow \frac{\bar{y}_2^2}{\bar{u}_2^2}\cdot \frac{\bar{y}_1^2}{\bar{u}_1^2}< 1 \Rightarrow \sqrt{\Big(\frac{\bar{y}_2}{\bar{u}_2}\cdot \frac{\bar{y}_1}{\bar{u}_1}\Big)^2}< 1\\[4pt]
   & \Rightarrow |H_2(0)H_1(0)|< 1\\[4pt]
   & \Rightarrow |\lambda_{max}(H_1(0)H_2(0))|<1.
\end{align*}
\end{remark}

\subsection{Closed-Loop Stability}
In the previous lemma we have constructed a candidate Lyapunov function of the positive feedback interconnection of the systems $H_1$ and $H_2$. We now introduce a robust stability result for this feedback interconnection as shown in Figure \ref{fig1}. By using Lyapunov stability theory and LaSalle's invariance principle, we show that the states $x_1$ and $x_2$ asymptotically approach the origin.

To establish this result, we represent the corresponding closed-loop system with the following state-space representation:
\begin{equation}\label{AUG}
   \dot{z}(t):=\varrho(z(t)), \quad  z(t):=\left[ \begin{matrix}
               x_{1}(t)\\[2pt]
                  x_{2}(t)
               \end{matrix}\right]\in \mathbb{R}^{2n}
\end{equation}
where the function $\varrho: \mathbb{R}^{2n}\rightarrow\mathbb{R}^{2n}$ is locally Lipschitz and $\varrho(0)=0$. Now we state the stability result in the following theorem.

\tikzstyle{block} = [draw, thick, rectangle,
    minimum height=4em, minimum width=4em]
\tikzstyle{sum} = [draw, circle, node distance=1cm]
\tikzstyle{input} = [coordinate]
\tikzstyle{output} = [coordinate]
\tikzstyle{pinstyle} = [pin edge={to-,thin,black}]
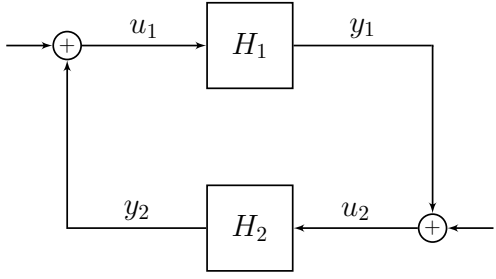
\begin{figure}[t]
\centering
\resizebox{6.5cm}{!}{%
\begin{tikzpicture}[auto, thick, node distance=2cm,>=latex']
    \node [input, name=input] {};
    \node [sum, right of=input, inner sep=1pt, minimum size=1pt,] (sum1) {$+$};
    \node [block, right of=sum1,
            node distance=3cm] (system) {\large $H_1$};
\node[inner sep=0,minimum size=0,right of=system,node distance=3cm] (k) {}; 
    \node [block, below of=system, node distance=3cm] (measurements) {\large $H_2$};
    \node[inner sep=0,minimum size=0,left of=measurements,node distance=3cm] (p) {}; 
    \node [sum, below of=k, inner sep=1pt, minimum size=1pt, node distance=3cm] (sum2) {$+$};
    \node [input,  right of=sum2, node distance=1cm](input2) {};
    \draw [draw,->] (input) -- node {\iffalse$r_1=0$\fi}node[pos=1.19]{\iffalse$+$\fi} (sum1);
    \draw [->] (sum1) -- node { \large $u_1$} (system);
    \draw [-] (system) -- node {\large $y_1$}(k);
    \draw [->] (k) -| node {} (sum2);
    \draw [->] (input2) -- node[above]{\iffalse\large $r_2$\fi} (sum2);
    \draw [->] (sum2) -- node[above]{\large $u_2$} (measurements);
    \draw [-] (measurements) -- node[above]{\large $y_2$} (p);
    \draw [->] (p) -|  node[pos=0.99]{\iffalse$-$\fi} (sum1);
\end{tikzpicture}}
\caption{Positive feedback interconnection of an NNI system $H_1$ and a WS-NNI system $H_2$.} \label{fig1}
\end{figure}
\begin{theorem}\label{main}
Consider the positive feedback interconnection of systems $H_1$ and $H_2$ as in Figure \ref{fig1}. Suppose that the system $H_1$ is  NNI and zero-state observable, and the system $H_2$ is WS-NNI. Moreover, suppose that Assumptions \ref{I}-\ref{IV} are satisfied. Then, the equilibrium point $z=0$ of the corresponding closed loop system \eqref{AUG} is asymptotically stable.
\end{theorem}
\begin{pf}
Consider the function $W(x_1,x_2)$ given in \eqref{lyap_fn} as a candidate Lyapunov function for the closed loop system \eqref{AUG}. Differentiating $W(x_1,x_2)$ with respect to $t$ and noting that $y_1=u_2$ and $y_2=u_1$, we get
\begin{equation*}
    \dot{W}(x_1,x_2)=\dot{V}_1(x_1)+ \dot{V}_2(x_2)-\dot{y}_1^T u_1 -\dot{y}_2^T u_2\leq0,
\end{equation*}
it follows that the closed-loop system \eqref{AUG} is at least Lyapunov stable.
\noindent Next, we use LaSalle's invariance principle to prove that the state trajectories of system \eqref{AUG} asymptotically approach to zero; that is, $(x_1,x_2)\rightarrow 0$. Since the system $H_2$ is WS-NNI, the system  $\dot{x}_2=f_2(x_2,0)$, is globally asymptotically stable (from Definition \ref{wsnni}). By using the result of \cite{eduardo1989remarks}, we have $\lim_{t\rightarrow\infty}x_2(t)=0$ which in turn leads to $\lim_{t\rightarrow\infty}y_2(t)=0$. For trajectories along which $\dot{W}(x_1,x_2)=0$, we have
\begin{equation}
   \dot{V}_2(x_2)-\dot{y}_2^T u_2=0.
\end{equation}
Thus, $\lim_{t\rightarrow\infty}u_2(t)=c$, since $H_2$ is WS-NNI.

\noindent Now, let $\Omega(z_0)$ be an $\omega$-limit set of a trajectory $z(t,z_0)$ with $\dot{W}(z)=0$. We show that $\Omega(z_0)=\{0\}$. For any $\alpha\in \Omega(z_0)$, we write
\begin{equation}
 \quad \alpha = \left[\begin{matrix}
                 \alpha_{1} \\[2pt]
                 \alpha_{2}
  \end{matrix}\right]\in \mathbb{R}^{2n}.
\end{equation}
Since $x_2(t)\rightarrow 0$, $\alpha_2=0$. The limit set  $\Omega(z_0)$ is an invariant set of the system \eqref{AUG}. In other words, $z(t,\alpha)\in \Omega(z_0)$ for all $t\geq0$. Then, $\alpha_2=0$ implies $z_2(t,\alpha)\equiv0$ and
\begin{equation}
  {\dot{z}}_1(t,\alpha) = f_1(z_1(t,\alpha),\,0), \quad 0\equiv h_1(z_1(t,\alpha)).
\end{equation}
\noindent Since the system $H_1$ is zero-state observable, we conclude that $z_1(t,\alpha)\equiv 0$. Thus, we see $\alpha=0$  and $\Omega(z_0)=\{0\}$. Hence, by Lemma 4.1 of \cite{khalil2002nonlinear}, $z(t,z_0)$ approaches $\Omega(z_0)=\{0\}$ as $t\rightarrow\infty$. We have $W(z(t,z_0))\equiv W(z_0)\rightarrow0$, as $t\rightarrow\infty$ and hence $W(z_0)=0$. It follows from the positive definiteness of $W$ that $z_0=0$. It is concluded that any bounded trajectory $z(t)$ satisfying $\dot{W}(z)\equiv0$ is the trajectory $z(t)\equiv0$. Now from the Lasalle invariance principle, any bounded trajectory $z(t)$ tends to the origin. Therefore, the origin $z=0$ is asymptotically stable.
\begin{flushright}
$\blacksquare$
\end{flushright}
\end{pf}
%
%

\section{Illustrative Example: Nonlinear Mass-Spring-Damper System}

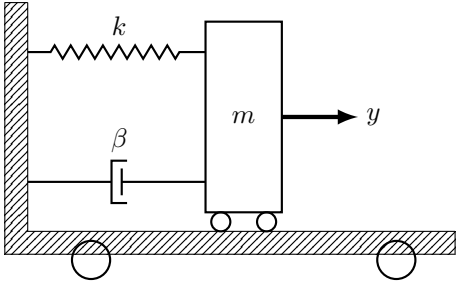
\begin{figure}[t!]
\centering
\usetikzlibrary{calc,patterns,decorations.pathmorphing,decorations.markings}
\resizebox{6cm}{!}{%
\begin{tikzpicture}
\tikzstyle{spring}=[thick,decorate,decoration={zigzag,pre length=0.3cm,post length=0.3cm,segment length=6}]
\tikzstyle{damper}=[thick,decoration={markings,
  mark connection node=dmp,
  mark=at position 0.5 with
  {
    \node (dmp) [thick,inner sep=0pt,transform shape,rotate=-90,minimum width=15pt,minimum height=3pt,draw=none] {};
    \draw [thick] ($(dmp.north east)+(2pt,0)$) -- (dmp.south east) -- (dmp.south west) -- ($(dmp.north west)+(2pt,0)$);
    \draw [thick] ($(dmp.north)+(0,-5pt)$) -- ($(dmp.north)+(0,5pt)$);
  }
}, decorate]
\tikzstyle{ground}=[fill,pattern=north east lines,draw=none,minimum width=0.75cm,minimum height=0.3cm,inner sep=0pt,outer sep=0pt]

\node [style={draw,outer sep=0pt,thick}] (M) [minimum width=1cm, minimum height=2.5cm] {$m$};

\node (ground) [ground,anchor=north,yshift=-0.25cm,minimum width=5.6cm,xshift=-0.03cm] at (M.south) {};
\draw (ground.north east) -- (ground.north west);
\draw (ground.south east) -- (ground.south west);
\draw (ground.north east) -- (ground.south east);

\node (fill) [ground,xshift=-0.15cm,minimum height = 0.3cm, minimum width = 0.3cm] at (ground.west) {};
\draw (fill.north west) -- (fill.south west);
\draw (fill.south west) -- (fill.south east);

\draw [thick] (M.south west) ++ (0.2cm,-0.125cm) circle (0.125cm)  (M.south east) ++ (-0.2cm,-0.125cm) circle (0.125cm);
\draw [thick] (M.south west) ++ (2.5cm,-0.625cm) circle (0.25cm)  (M.south east) ++ (-2.5cm,-0.625cm) circle (0.25cm);
\node (wall) [ground, rotate=-90, minimum width=3cm,anchor=south east] at (fill.north west) {};
\draw (wall.north east) -- (wall.north west);
\draw (wall.north west) -- (wall.south west);
\draw (wall.south west) -- (wall.south east);

\node (y) at (M.east) [xshift = 1.2cm] {$y$};

\draw [spring] (wall.170) -- ($(M.north west)!(wall.170)!(M.south west)$);
\draw [damper] (wall.10) -- ($(M.north west)!(wall.10)!(M.south west)$);
\node (b) at (wall.10) [xshift = 1.2cm,yshift=0.55cm] {$\beta$};
\node (k) at (wall.170) [xshift = 1.2cm,yshift=0.35cm] {$k$};

\draw [-latex,ultra thick] (M.east) ++ (0cm,0cm) -- +(1cm,0cm);
\end{tikzpicture}}
\caption{Nonlinear Mass Spring Damper System.}\label{msd}
\end{figure}

\noindent To illustrate the applicability of the above NNI stability result, we consider the problem of control design for a nonlinear mass-spring-damper system as shown in Figure \ref{msd}, which has been shown to be NNI. We aim at designing a controller by virtue of Theorem \ref{main} to robustly stabilize the system. 
Here, the system is assumed to be nonlinear and obey the force law
\begin{equation}
 f=2(x+x^3),
\end{equation}
where $x$ is the displacement of the spring. From Newton's second law, the system is governed by the equation
\begin{equation}\label{mass_spring}
  m\ddot{x}+\beta \dot{x}+k(x+x^3)=f(t)=u(t),
\end{equation}
Setting $\beta=1,~ m=1,~ k=2$ and define states $x_1=x, \ x_2=\dot{x}$, we obtain the following nonlinear state-space representation
\begin{align}\label{matrixform}
\dot{x}&=\left[ \begin{matrix}
                  \dot{x}_{1} \\
                  \dot{x}_{2} \\
               \end{matrix}\right]=
               \left[ \begin{matrix}
                  x_{2} \\
                  \frac{-k}{m}(x_1+x_1^3)-\frac{\beta}{m}x_2+\frac{u(t)}{m}\\
               \end{matrix}\right];\\
               y&= \left[\begin{matrix}
                1 \hspace{1.5mm} 0
               \end{matrix}\right]
                \left[ \begin{matrix}
                  x_{1} \\
                  x_{2} \\
               \end{matrix}\right].
\end{align}
\begin{figure}[t!]
\centering
\tikzstyle{block} = [draw, thick, rectangle,
    minimum height=2em, minimum width=4em]
\tikzstyle{sum} = [draw, circle,inner sep=0pt,minimum size=1pt, node distance=1cm]
\tikzstyle{input} = [coordinate]
\tikzstyle{output} = [coordinate]
\tikzstyle{pinstyle} = [pin edge={to-,thin,black}]
\tikzstyle{int}=[draw]
\resizebox{7.3cm}{!}{%
\begin{tikzpicture}[auto,>=latex']  
\tikzstyle{block} = [draw, thick, rectangle,
    minimum height=3em, minimum width=4em, text width=2cm]
    \node [int] (a)[align=left] {
     \tiny{\hspace{.24em} MSD} \\ \tiny{System}};
    \node (b) [left of=a,node distance=1.2cm, coordinate] {a};
    \node [int] (c) [right of=a,  node distance=2.8cm, align=left] {
     \tiny{\hspace{.88em} SNI} \\ \tiny{Controller}};
    \node [coordinate] (end) [right of=c, node distance=1.4cm]{};
    \path[->] (b) edge node {\small$\bar{u}$} (a);
    \path[->] (a) edge node {\small$\bar{y}  \hspace{.65cm} \bar{u}_c$} (c);
    \draw[->] (c) edge node {\small$\bar{y}_c$} (end) ;
\end{tikzpicture}}
\caption{Open-loop system of the MSD and the SNI controller \eqref{sni_cont}  ('c' refers to the controller) in the steady-state case.}\label{sys3}
\end{figure}
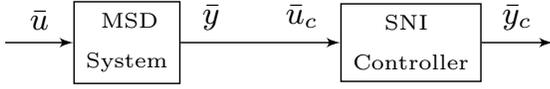
\noindent Now suppose the system is controlled with the SISO integral resonant controller
\begin{equation}\label{sni_cont}
  C(s)= \frac{\Gamma}{s+\Gamma\Phi}\CommaBin
\end{equation}
where $\Gamma, \Phi$ are positive constants. The transfer function $C(s)$ is strictly negative imaginary \cite{petersen2010}.  We first check the assumptions of the open-loop interconnection required by Theorem \ref{main}. From the state equation \eqref{matrixform}, when $\dot{x}=0$, we get
\begin{align}
  \bar{x}_2=0; \\
\label{Cubiceq} \bar{x}_1^3+\bar{x}-\frac{\bar{u}}{2}=0;
\end{align}
Using the discriminant method to solve the cubic equation, since we have
\begin{equation}
\Delta=-4-27\Big(\frac{\bar{u}}{k}\Big)^2<0,
\end{equation}
and therefore equation \eqref{Cubiceq} has one real solution. Using Cardano's cubic formula, we get the solution
\begin{equation*}
\bar{y}=\bar{x}=\displaystyle\sqrt[3]{\frac{\bar{u}}{2k}+\sqrt{\frac{1}{27}+\frac{\bar{u}^2}{4k^2}}}+\sqrt[3]{\frac{\bar{u}}{2k}-\sqrt{\frac{1}{27}+\frac{\bar{u}^2}{4k^2}}}
\end{equation*}
and from the open-loop interconnection, in Figure~\ref{sys3}, we have
$$
\bar{y}_c=\frac{1}{\Phi}\bar{u}_c=\frac{1}{\Phi}\bar{y},
$$
thus we obtain
\begin{equation*}
\bar{y}_c=\frac{1}{\Phi}\cdot\left(\sqrt[3]{\frac{\bar{u}}{2k}+\sqrt{\frac{1}{27}+\frac{\bar{u}^2}{4k^2}}}+\sqrt[3]{\frac{\bar{u}}{2k}-\sqrt{\frac{1}{27}+\frac{\bar{u}^2}{4k^2}}}\ \right).
\end{equation*}
\noindent Plotting this function, we can choose a value for the controller parameter $\Phi$ such that the sector bound condition
\begin{equation}\label{sec1}
  \bar{y}_c^2\leq \gamma^2 \bar{u}^2
\end{equation}
is satisfied for $0<\gamma<1$, as shown in Figure \ref{fff}.
%
%
\begin{figure}[H]
  \centering
  \includegraphics[scale=0.315]{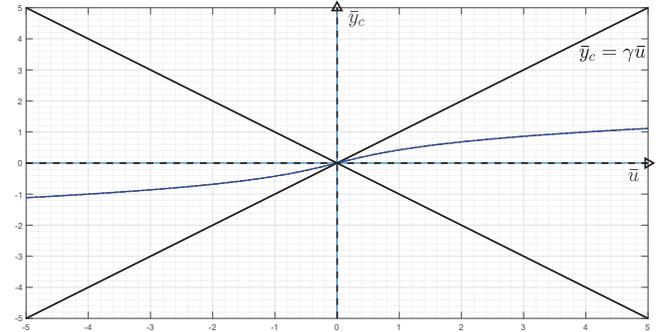}
  \caption{The sector bound condition \eqref{sec1} is satisfied with $\Phi=2$ for NNI controller design of the MSD system. \label{fff}}
\end{figure}

\begin{remark}
In recent research work, in view of Theorem \ref{main}, a nonlinear negative imaginary approach has been used to design velocity-free controller for the problems quadrotors attitude control \cite{ghallab2021velocity}, quadrotors position control \cite{ghallab2021quadrotors}, and robotic manipulators control problem \cite{ghallab2021manipulators}. In these papers, the NNI control design approach aims at replacing rate sensors on the quadrotor\slash manipulator structure with displacement sensors.
\end{remark}

\section{Nonlinear Negative Imaginary System with Free Motion}
In this section, we aim to further generalize the notion of nonlinear negative imaginary to establish the nonlinear analog of negative imaginary systems with free body motion presented in \cite{Mabrok2013} for LTI systems. In \cite{Mabrok2013}, the definition and the stability result for LTI systems was extended to the case of NI systems with a pure integrator (i.e. a pole at the origin). Motivated by the results of \cite{Mabrok2013}, we seek to extend Definition \ref{nonlinear_NI} to establish the conditions under which a nonlinear systems in a cascade connection with an integrator would be nonlinear negative imaginary.

\subsection{Nonlinear Negative System with an Integrator}

We consider here a SISO nonlinear system which is affine in the input of the form
\begin{equation}\label{nl1}
\Sigma :
\left\{
  \begin{array}{ll}
    \dot{\eta} = f(\eta)+g(\eta)\xi \\
   \dot{\xi} = u \\
    y= h(\eta)
  \end{array}
\right.
\end{equation}
where $\eta \in \mathbb{R}^n$, $u\in\mathbb{R}^n$ is the scalar control input, and the functions $f:\mathbb{R}^n\rightarrow\mathbb{R}^n, \ g:\mathbb{R}^n\rightarrow\mathbb{R}^n, \  h:\mathbb{R}^n\rightarrow\mathbb{R}$, are continuously differentiable, where $f(0)=0$, and $g(\eta)\neq0$ for all  $\eta$. This system is considered as a cascade connection of an integrator with the nonlinear $\eta$ subsystem as shown in Figure \ref{single_integrator}.

\noindent The idea now is to try to find storage function which will satisfy the dissipation inequality \eqref{lyapD} according to Definition \ref{nonlinear_NI}.
\noindent That is, we need to ensure that there exists a continuously differentiable storage function $V(\eta,\xi)$ which is positive-definite and its time-derivative along trajectories of the system \eqref{nl1} satisfies the following inequality
\begin{align*}
  \dot{V}(\eta,\xi) &\leq \dot{y}u= \nabla h(\eta)[f(\eta)+g(\eta)\xi]u \hspace{.45cm} \forall \ u,  \xi,  \eta.
\end{align*}
Now, for all $u, \xi, \eta$, we obtain
\begin{align*}
  \frac{\partial V}{\partial\eta}(f(\eta)+g(\eta)\xi)+ \frac{\partial V}{\partial\xi}u
     \leq u\nabla h(\eta) f(\eta)+u\nabla h(\eta)g(\eta)\xi.
\end{align*}
Rearranging the terms of the above inequality, we have, for all $\ u, \xi, \eta$, the following inequality
\begin{equation*}
 \frac{\partial V}{\partial\eta}[f(\eta)+g(\eta)\xi]
  +u[\frac{\partial V}{\partial\xi}-\nabla h(\eta) f(\eta)
   -\nabla h(\eta)g(\eta)\xi]\leq0.
\end{equation*}
Then we can say, $\forall \ t \in [0,\infty)$, the function $V$ satisfies $\dot{V} \leq \dot{y}(t)u(t)$ if and only if the following two conditions are satisfied for all $\xi, \eta$:
\begin{align}
\label{rr1}&\frac{\partial V}{\partial\xi}=\nabla h(\eta)[f(\eta)+g(\eta)\xi];\\[4pt]
\label{rr3}& \hspace{.5cm} \frac{\partial V}{\partial\eta}[f(\eta)+g(\eta)\xi]\leq0.
\end{align}
\noindent Integrating both sides of \eqref{rr1} with respect to $\xi$, one has
\begin{equation*}\label{rr2}
 V(\eta,\xi)=\bar{V}(\eta)+\nabla h(\eta)f(\eta)\xi+\frac{1}{2}\nabla h(\eta)g(\eta)\xi^2\CommaBin
\end{equation*}
which is a quadratic function in $\xi$, and  the term $ \bar{V}(\eta)$ is a constant of integration which can be freely chosen such that
\begin{equation*}
  \bar{V}(\eta)\geq\frac{(\nabla h(\eta) f(\eta))^2}{2\nabla h(\eta)g(\eta)}\CommaBin
\end{equation*}
for all $\eta$, so we can choose
\begin{equation*}
\bar{V}(\eta)=\frac{(\nabla h(\eta) f(\eta))^2}{2\nabla h(\eta)g(\eta)}\mathbin{\raisebox{0.5ex}{.}}
\end{equation*}
\tikzstyle{int}=[draw, thick, rectangle, minimum height = 3em,
    minimum width = 3em]
\tikzstyle{init} = [pin edge={to-,thick,black}]
\begin{figure}[t!]
\centering
\setlength{\belowcaptionskip}{0pt}
\resizebox{8.5cm}{!}{
\begin{tikzpicture}[node distance=2.5cm,auto,>=latex']
    \node [int] (a) {$\displaystyle \int$};
    \node (b) [left of=a,node distance=2cm, coordinate] {a};
    \node [int] (c) [right of=a,  node distance=3cm] {$\begin{matrix}
    \dot{\eta} = f(\eta)+g(\eta)\xi \\
    y= h(\eta)
  \end{matrix}$};
    \node [coordinate] (end) [right of=c, node distance=3cm]{};
    \path[->] (b) edge node {$u$} (a);
    \path[->] (a) edge node {$\xi$} (c);
    \draw[->] (c) edge node {$y$} (end) ;
\end{tikzpicture}}
\caption{Cascade interconnection of an input-affine nonlinear system with an integrator.} \label{single_integrator}
\end{figure}
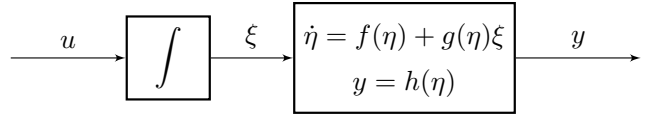
Now, we find conditions on the functions $f, g$ and $h$ for which the conditions \eqref{rr1}, \eqref{rr3} are satisfied and ensure at the same time that the function $V(\eta, \xi)$ is positive-definite, and hence the system \eqref{nl1} is NNI according to Definition \ref{nonlinear_NI}. This is made clear in the next theorem.

\begin{theorem}\label{singINT}
Consider a nonlinear system, affine in the input, of the form \eqref{nl1}. 
Assume that the following conditions are satisfied for ll $\eta$:
\begin{itemize}
\item [\emph{(1)}] $\nabla h(\eta) g(\eta)=constant>0$;\\[2pt] 
\item [\emph{(2)}] $(\nabla h) \nabla (\nabla h(\eta) f(\eta))\leq0;$ \\[2pt]
\item [\emph{(3)}] $\nabla h(\eta) f(\eta)\neq 0, \quad \forall \ \eta \neq 0$.

\end{itemize}
Then the system \eqref{nl1} is nonlinear negative imaginary with a positive-definite storage function $V(\eta,\xi)$ given by
\begin{equation}\label{storage2}
  V(\eta,\xi) = \frac{(\nabla h(\eta) f(\eta))^2}{2\nabla h(\eta) g(\eta)}+\nabla h(\eta)f(\eta)\xi+\frac{1}{2}\nabla h(\eta)g(\eta)\xi^2.
\end{equation}
\end{theorem}
\begin{pf}
From \eqref{storage2}, it can be easily seen that $V(0,0)=0$. Now, for $(\eta,\xi)\neq0$ we can see that the storage function $V(\eta,\xi)$ is a quadratic function in $\xi$ with its discriminant $D$ given by:
\begin{equation*}
  D=(\nabla h(\eta)f(\eta))^2-\frac{(\nabla h(\eta) f(\eta))^2}{\nabla h(\eta) g(\eta)}\nabla h(\eta)g(\eta)=0
\end{equation*}
and from condition $(1)$, we have $\nabla h(\eta) g(\eta)>0$ which implies that $V(\eta,\xi)\geq 0$ for all $(\eta, \xi)\neq0$. To ensure positive definiteness of the function $V(\eta,\xi)$, we rewrite it in the following way:
\begin{align*}
  V(\eta, \xi)&=\frac{(\nabla h(\eta) f(\eta))^2}{2\nabla h(\eta) g(\eta)}+\nabla h(\eta)f(\eta)\xi+\frac{1}{2}\nabla h(\eta)g(\eta)\xi^2\\[6pt]
  & =\frac{1}{2}\Big(\frac{(\nabla h(\eta) f(\eta))}{\sqrt{2}\sqrt{\nabla h(\eta) g(\eta)}}+\frac{1}{\sqrt{2}}\sqrt{\nabla h(\eta)g(\eta)}\xi\Big)^2\\[4pt]
 & +\frac{1}{4}\frac{(\nabla h(\eta) f(\eta))^2}{\nabla h(\eta) g(\eta)} + \frac{1}{4}\nabla h(\eta)f(\eta)\xi^2,
\end{align*}
and from condition $(4)$, we can see that $V(\eta, \xi)\neq 0$ for all $(\eta, \xi)\neq0$; that is, the function  $V(\eta, \xi)$ is positive-definite.

Next, From \eqref{storage2}, differentiating the function $V(\eta,\xi)$ with respect to $t$, we get
\begin{align*}
\dot{V}(\eta,\xi)&=\frac{\partial V}{\partial\eta}\dot{\eta}+ \frac{\partial V}{\partial\xi}\dot{\xi}\\
  &=\frac{\nabla(\nabla h(\eta) f(\eta))^2}{2\nabla h(\eta)
  g(\eta)}\dot{\eta}+\nabla(\nabla h(\eta)f(\eta))\xi \dot{\eta}\\
  &+\frac{1}{2}\nabla(\nabla h(\eta)g(\eta))\xi^2\dot{\eta} + \nabla h(\eta) f(\eta)\dot{\xi}\\[2pt] & \quad +\nabla h(\eta) g(\eta)\xi\dot{\xi}\\[2pt]
  &= \alpha(f(\eta)+ g(\eta)\xi)+\beta (f(\eta)+ g(\eta)\xi)\xi\\[2pt]
  & \quad + \nabla h(\eta) (f(\eta)+ g(\eta)\xi) u
\end{align*}
where
\begin{equation}
\alpha = \frac{\nabla(\nabla h(\eta) f(\eta))^2}{2\nabla h(\eta) g(\eta)}\CommaBin \quad \beta = \nabla (\nabla h(\eta) f(\eta)).
\end{equation}
Hence, we get
\begin{equation*}
\dot{V}(\eta,\xi)= \alpha f(\eta)+(\beta f(\eta)+\alpha g(\eta))\xi+\beta g(\eta)\xi^2+\dot{y}(t)u(t)
\end{equation*}
The quantity $\alpha f(\eta)+(\beta f(\eta)+\alpha g(\eta))\xi+\beta g(\eta)\xi^2$ is a quadratic function in $\xi$, and from conditions $(2)$, $(3)$, it is negative semi-definite. Therefore, we have
\begin{equation}
  \dot{V}(\eta,\xi)\leq \dot{y}(t)u(t)
\end{equation}
for all $ t\geq 0$, which shows that the system \eqref{nl1} is nonlinear negative imaginary. Then the proof completes.
\begin{flushright}
$\blacksquare$
\end{flushright}
\end{pf}

To illustrate the result of Theorem \ref{singINT}, we consider the following examples.

\begin{example} Consider the following nonlinear system
\begin{equation}\label{nl220}
\Sigma :
\left\{
  \begin{array}{ll}
    \dot{\eta}_1 = \eta_2; \\
    \dot{\eta}_2 = -\eta_1^3-\eta_2+\xi;\\
   \dot{\xi} = u; \\
    y= \eta_2,
  \end{array}
\right.
\end{equation}
which has the form of system \eqref{nl1} where
\begin{equation*}
f(\eta)=\left[
\begin{array}{c}
\eta_2\\
-\eta_1^3-\eta_2\\
\end{array}
\right], \quad
 g(\eta)=\left[
\begin{array}{c}
0\\
1\\
\end{array}
\right],
\end{equation*}
and
\begin{equation*}
h(\eta)= \left[
\begin{array}{r}
0\hspace{1.5mm} 1
\end{array}
\right]\left[
\begin{array}{c}
\eta_1\\
\eta_2\\
\end{array}
\right].
\end{equation*}
We can see that
$$
\nabla h(\eta) g(\eta)= \left[
\begin{array}{r}
0\hspace{1.5mm} 1
\end{array}
\right]\left[
\begin{array}{c}
0\\
1\\
\end{array}
\right]=1>0,
$$

$$
\nabla[\nabla h(\eta) g(\eta)]g(\eta)=0,
$$
and also
$$
\nabla h \nabla (\nabla h(\eta) f(\eta))= \left[
\begin{array}{r}
0\hspace{1.5mm} 1
\end{array}
\right]\left[
\begin{array}{c}
-3\eta_1^2\\
-1\\
\end{array}
\right]=-1\leq0.
$$
In virtue of Theorem \ref{singINT}, we conclude that the system \eqref{nl220} is nonlinear negative imaginary.
\begin{flushright}
$\blacksquare$
\end{flushright}
\end{example}

As a special result of Theorem \ref{singINT}, we can show that a cascade connection of a positive-real system and an integrator will result in NI systems, which in fact comply with the result of Lemma \ref{pr_ni}. We present the following two examples.

%

\begin{example}
Consider the following LTI system
\begin{equation}\label{nl01}
\left\{
\begin{array}{ll}
\dot{\eta}(t) = a~\eta(t)+b~\xi; \\
         \hspace{.42cm} \dot{\xi} = u;\\
          y(t)= c~\eta(t),
\end{array}
\right.
\end{equation}
where $\eta \in \mathbb{R}^n$, $u \in \mathbb{R}$, $ y\in \mathbb{R}$, and $a$, $b$, $c$ are real constants. The $\eta$ subsystem is assumed to be positive real with a transfer function $G(s)$ given by
\begin{equation}\label{pr1}
  G(s)=\frac{cb}{s-a}.
\end{equation}
where $c\ b>0$ and $a<0$. It can be easily verify that the assumptions of Theorem \ref{singINT} are satisfied. Then, we have the following positive-definite storage function
\begin{align*}
V(\eta,\xi)&=\frac{1}{2}\frac{\nabla h(\eta) f(\eta)^2}{g(\eta)}+\nabla h(\eta)f(\eta)\xi+\frac{1}{2}\nabla h(\eta)g(\eta)\xi^2\\
             &=\frac{1}{2}\frac{ca^2}{b}\eta^2+ca\eta\xi+\frac{1}{2}cb\xi^2,
\end{align*}
which has a time-derivative evaluated as
\begin{align*}
   \dot{V}(\eta,\xi)&\leq \frac{a^2c}{b}\eta\dot{\eta} + ca\dot{\eta}\xi + ca\eta\dot{\xi} + cb\xi\dot{\xi}\\
&=\frac{ac}{b}(a\eta + b\xi)\dot{\eta} +c(a\eta + b\xi)\dot{\xi}\\
&=\frac{ac}{b}\dot{\eta}^2+c\dot{\eta}\dot{\xi}\leq c\dot{\eta}\dot{\xi}.
\end{align*}
Hence, that the system \eqref{nl01} is nonlinear negative imaginary.
\begin{flushright}
$\blacksquare$
\end{flushright}
\end{example}

\begin{example}
Consider the following positive real system
  \begin{equation}\label{pr}
    G(s)=\frac{s+1}{s^2+s+1}
  \end{equation}
which has the following state-space representation
\begin{equation}\label{sp}
A=\left[
\begin{array}{cc}
-1 & -1\\
1 & 0\\
\end{array}
\right], \
B=\left[
\begin{array}{c}
1\\
0\\
\end{array}
\right], \
C=\left[
\begin{array}{r}
1\hspace{1.5mm} 1
\end{array}
\right], \
D=0.
\end{equation}
This state-space realization can be put in the form of the system \eqref{nl1} where
\begin{equation}
  f(\eta)= A\eta, \quad h(\eta)=C\eta, \quad  g(\eta)=B.
\end{equation}
We can see that,
$$\nabla h(\eta) g(\eta)=CB=1>0,$$
$$\nabla[\nabla h(\eta) g(\eta)]g(\eta)=0,$$ 
and
$$\nabla h \nabla (\nabla h(\eta) f(\eta))=\left[
\begin{array}{r}
1\hspace{1.5mm} 1
\end{array}
\right]\left[
\begin{array}{c}
0\\
-1\\
\end{array}
\right]=-1\leq0,$$
which shows that the Assumption of Theorem \ref{singINT} are satisfied, and hence the system \eqref{pr} is shown to be a nonlinear negative imaginary.
\begin{flushright}
$\blacksquare$
\end{flushright}
\end{example}

\subsection{Robust Stability Result}
We conclude this section by introducing a stability robustness result a positive feedback interconnection of two NNI systems as shown in Figure~\ref{fig22}, where the system $\Sigma_1$ is of the form \eqref{nl1} and is given by,  
\begin{equation}\label{H_1}
 \displaystyle \Sigma_1{:} \ \left\{
  \begin{array}{ll}
     \dot{\eta}_1&=f_1(\eta_1)+g_1(\eta_1)\xi; \\
        \dot{\xi}&=u_1;\\
              y_1&=h_1(\eta_1);
         \end{array}
\right.
\end{equation}
\noindent and an input-affine nonlinear system $\Sigma_2$ given by,
\begin{equation}
\displaystyle \Sigma_2{:} \
\left\{
  \begin{array}{ll}
    \dot{\eta}_2&=f_2(\eta_2)+g_2(\eta_2)u_2;\\
     y_2&=h_2(\eta_2).
  \end{array}
\right.
\end{equation}
Here all the functions involved are assumed to be sufficiently smooth to grantees the existence of the solutions. The stability robustness of the feedback interconnection of systems $\Sigma_1$ and $\Sigma_2$ can be established in a similar way to that of Theorem \ref{main}.

In order to state the stability result for the feedback interconnection of $\Sigma_1$ and $\Sigma_2$, we represent the interconnection by the following state-space representation:
\begin{equation}\label{AUG1}
   \dot{z}(t):=\varrho(z(t)), \quad  z(t):=\left[ \begin{matrix}
               \eta_{1}(t)\\[2pt]
               \xi(t)\\[2pt]
                  \eta_{2}(t)
               \end{matrix}\right]\in \mathbb{R}^{2n+1},
\end{equation}
where $\varrho: \mathbb{R}^{2n+1}\rightarrow\mathbb{R}^{2n+1}$ is assumed to be locally Lipschitz and $\varrho(0)=0$. 

\tikzstyle{block} = [draw, thick, rectangle,
    minimum height=4em, minimum width=4em]
\tikzstyle{sum} = [draw, circle,inner sep=0pt,minimum size=1pt, node distance=1cm]
\tikzstyle{input} = [coordinate]
\tikzstyle{output} = [coordinate]
\tikzstyle{pinstyle} = [pin edge={to-,thin,black}]
\tikzstyle{my help lines}=[gray,thick,dashed]
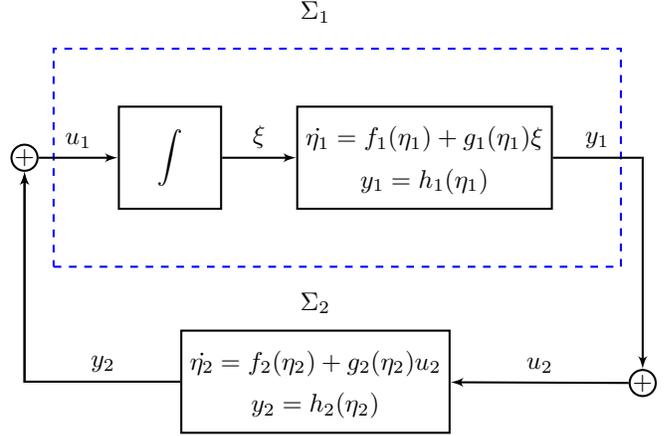
\begin{figure}[t!]
\centering
\setlength{\belowcaptionskip}{-8pt}
\resizebox{8.8cm}{!}{
\begin{tikzpicture}[auto, thick, node distance=2cm,>=latex']
    \node [input, name=input] {};
    \node [sum, right of=input] (sum1) {$+$};
   \node [block, right of=sum1, node distance=2cm] (integrator) {$\displaystyle \int$};
    \node [block, right of=integrator,
            node distance=3.5cm] (system) {$\begin{matrix}
    \dot{\eta_1} = f_1(\eta_1)+g_1(\eta_1)\xi \\
    y_1= h_1(\eta_1)
  \end{matrix}$};
  \draw[style=my help lines, blue] (1.4,-1.5)rectangle (9.2,1.5) ;
  \node[inner sep=0,minimum size=0](M1) at (5,2) {$\Sigma_1$};
  \draw [->] (integrator) -- node[name=u] {$\xi$} (system);
\node[inner sep=0,minimum size=0,right of=system,node distance=3cm] (k) {}; 
    \node [block, below of=M1, node distance=5.082cm] (measurements) {$\begin{matrix}
    \dot{\eta_2} = f_2(\eta_2)+g_2(\eta_2)u_2 \\
    y_2= h_2(\eta_2)
  \end{matrix}$};
   \node[inner sep=0,minimum size=0](M2) at (5,-2) {$\Sigma_2$};
    \node[inner sep=0,minimum size=0,left of=measurements,node distance=4cm] (p) {}; 
    \node [sum, below of=k, node distance=3.1cm] (sum2) {$+$};
    \node [input,  right of=sum2, node distance=1cm](input2) {};
    \draw [->] (sum1) -- node {$u_1$} (integrator);
    \draw [-] (system) -- node {$y_1$}(k);
    \draw [->] (k) -| node {} (sum2);
    \draw [->] (sum2) -- node[above]{$u_2$} (measurements);
    \draw [-] (measurements) -- node[above]{$y_2$} (p);
    \draw [->] (p) -|  node {} (sum1);
    \draw [->] (measurements) -| node  [pos=0.99] {}
     node [near end] {} (sum1);
\end{tikzpicture}}
\caption{Positive feedback interconnection of NNI systems $\Sigma_1$ and $\Sigma_2$.} \label{fig22}
\end{figure}

We have now the following robust stability result.
\begin{theorem}\label{nonlin_stab2}
Consider the positive feedback interconnection of the systems $\Sigma_1$ and $\Sigma_2$ as in Figure \ref{fig22}. Suppose that the system $\Sigma_1$ is NNI and zero-state observable, and the system $\Sigma_2$ is WS-NNI. Moreover, suppose that the Assumptions \ref{I}-\ref{IV} are satisfied. Then, the equilibrium point $z=0$ of the corresponding closed-loop system \eqref{AUG1} is asymptotically stable.
\end{theorem}
\begin{pf}
  The proof of this theorem follows by a similar argument as in the proof of Theorem \ref{main}.
  \begin{flushright}
$\blacksquare$
\end{flushright}
\end{pf}

\section{Conclusion}
In this paper, the notion of the negative imaginary systems has been generalized from the linear case to the nonlinear case. A time-domain definition of nonlinear negative imaginary systems has been introduced as a dissipative system property with an appropriate supply rate. Next, the stability robustness of a positive feedback interconnection where the plant corresponds to a nonlinear negative imaginary system and the controller corresponds to a weak strict negative imaginary system has been established under a set of technical assumptions.  A standard mechanical system involving a nonlinear mass-spring-damper system has been considered to illustrate this nonlinear NI stability result. Furthermore, the notion of nonlinear negative imaginary systems has been further extended to include the case of free motion. A cascade connection of a nonlinear system, affine in the input, and integrator has been shown to be nonlinear NI under a set of assumptions. Finally, we introduced a stability result for a positive feedback system where the plant is nonlinear NI with an integrator and the controller is a weakly strictly negative imaginary system.


\bibliographystyle{plain}        
\bibliography{ifacbib2}           



\appendix
\section{Appendix} \label{append}  
We show for an LTI system of the form \eqref{eq:xdot1}-\eqref{eq:y1}, the WS-NNI property reduces to the SNI property of the system. This can be readily seen by considering a storage function $V(x)=\frac{1}{2}x^TPx$, where $x$ is the state vector of the system and the matrix $P$ is a positive-definite symmetric matrix which satisfies the LMI \eqref{rr}. Differentiating the function $V$ with respect to the time $t$, one has
\begin{align*}
\dot{V}&(x(t))\\ &= \frac{1}{2}x^TP\dot{x}+\frac{1}{2}\dot{x}^TPx\\
                &= \frac{1}{2}x^TP(Ax+Bu)+\frac{1}{2}(Ax+Bu)^TPx\\
                &= \frac{1}{2}x^T(PA+A^TP)x+x^TPBu\\
                &= -\frac{1}{2}x^T L^TLx+u^TB^TPx\\
                &= -\frac{1}{2}x^T L^TLx+u^T(CA-W^TL)x\\
                &= -\frac{1}{2}x^T L^TLx+u^T(CAx+CBu)-u^TCBu-u^TW^TLx\\
                &=-\frac{1}{2}x^TL^TLx+u^T\dot{y}-u^T(CB+CB^T)u-u^TW^TLx\\
                &= -\frac{1}{2}x^T L^TLx+u^T\dot{y}-u^TWW^Tu-u^TW^TLx\\
                &=\dot{y}^Tu-(LPx+LC^{T}u)^{T}(LPx+LC^{T}u)\\
                &=\dot{y}^Tu-\tilde{y}^{T}\tilde{y}
\end{align*}
where $\tilde{y}$ is the output of the auxiliary system
\begin{align}
\dot{x}&= Ax+Bu\\
\tilde{y}&=LPx+LC^{T}u
\end{align}
whose transfer function is 
$$
W(s)=\begin{bmatrix}
\begin{array}{c|c}
A & B \\ \hline LP \ & LPA^{-1}B
\end{array}
\end{bmatrix}.
$$
We have $\tilde{Y}(s)=W(s)U(s)$ where according to Lemma \ref{SNI}, the matrix $W(s)$ has no zeros on the imaginary axis except possibly at the origin. Whenever $\dot{V}=\dot{y}^Tu$, then $\tilde{y}\equiv0$, and thus $u(t)$ can only either exponentially diverge or exponentially converge to zero. However, if $u(t)$ diverges, it follows from $\tilde{y}(t)=LPx+LC^{T}u \equiv 0$ that $x$ should also diverge which is a contradiction. It then follows that $u$ should converge to zero, which proves that the LTI system \eqref{eq:xdot1}-\eqref{eq:y1} is WS-NNI according to Definition \ref{wsnni}.

\section{Acknowledgment}
The authors would like to thank Kanghong Shi, Australian National University, for helpful discussions and suggestions in the context of Ghallab et al. (2018).
\end{document}